\documentclass[12pt]{amsart}
\usepackage{amsmath,amssymb,latexsym,amsthm,amscd}

\theoremstyle{plain}
\newtheorem{theorem}[subsection]{Theorem}
\newtheorem{proposition}[subsection]{Proposition}
\newtheorem{corollary}[subsection]{Corollary}
\newtheorem{lemma}[subsection]{Lemma}
\newtheorem{conjecture}[subsection]{Conjecture}

\theoremstyle{definition}

\theoremstyle{remark}
\newtheorem{remark}[subsection]{Remark}

\newcommand{\Ebar}{{\overline{E}}}
\newcommand{\Lambdabar}{{\overline{\Lambda}}}

\newcommand{\bC}{{\mathbb{C}}}
\newcommand{\bF}{{\mathbb{F}}}
\newcommand{\bFbar}{{\overline{\mathbb{F}}}}
\newcommand{\bQ}{{\mathbb{Q}}}
\newcommand{\bQbar}{{\overline{\mathbb{Q}}}}
\newcommand{\bQell}{{\mathbb{Q}_\ell}}
\newcommand{\bQellbar}{{\overline{\mathbb{Q}}_\ell}}
\newcommand{\bR}{{\mathbb{R}}}
\newcommand{\bZ}{{\mathbb{Z}}}

\newcommand{\Fbold}{{\mathbf{F}}}
\newcommand{\Lbold}{{\mathbf{L}}}

\newcommand{\tD}{{\operatorname{D}}}
\newcommand{\tH}{{\operatorname{H}}}
\newcommand{\tK}{{\operatorname{K}}}
\newcommand{\tR}{{\operatorname{R}}}
\newcommand{\tX}{{\operatorname{X}}}

\newcommand{\ch}{{\text{\rm ch}}}

\newcommand{\Cht}{{\text{Cht}^{r,\overline{p}\le p}_N}}
\newcommand{\Chtbar}{{\overline{\Cht}}}
\newcommand{\Xfrak}{{\mathfrak{X}}}

\newcommand{\etabar}{{\bar{\eta}}}
\newcommand{\xbar}{{\bar{x}}}
\newcommand{\kbar}{{\overline{k}}}
\newcommand{\kx}{{\kappa(x)}}
\newcommand{\kxbar}{{\kappa(\xbar)}}

\newcommand{\Utilde}{{\widetilde{U}}}
\newcommand{\Vtilde}{{\widetilde{V}}}
\newcommand{\ptilde}{{\widetilde{p}}}

\newcommand{\Fcal}{{\mathcal{F}}}
\newcommand{\Gcal}{{\mathcal{G}}}
\newcommand{\Hcal}{{\mathcal{H}}}
\newcommand{\Lcal}{{\mathcal{L}}}

\newcommand{\Garith}{{\operatorname{G_{arith}}}}
\newcommand{\Ggeom}{{\operatorname{G_{geom}}}}

\newcommand{\GL}{{\operatorname{GL}}}
\newcommand{\Mult}{{\operatorname{Mult}}}

\newcommand{\Der}{{\operatorname{Der}}}
\newcommand{\Fr}{{\operatorname{Fr}}}
\newcommand{\Frob}{{\operatorname{Frob}}}
\newcommand{\Gal}{{\operatorname{Gal}}}
\newcommand{\Irr}{{\operatorname{Irr}}}
\newcommand{\Spec}{{\operatorname{Spec}}}

\newcommand{\cris}{{\text{cris}}}
\newcommand{\gr}{{\text{gr}}}
\newcommand{\semisimp}{{\text{ss}}}

\newcommand{\larrow}{{\kern4pt\leftarrow\kern-12pt}}
\newcommand{\rarrow}{{\kern-12pt\rightarrow\kern4pt}}

\newcommand{\hookuparrow}{{\begin{matrix} \uparrow \\[-.6ex]
				\kern-.31em{}^{\cup} \\[-1.35ex] \end{matrix} }}

\newcommand{\hookdownarrow}{{\begin{matrix} \kern-.31em{}^{\cap} \\[-1.95ex]
				\downarrow \\ \end{matrix} }}

\begin{document}

\title[Independence of~${\ell}$ of Monodromy Groups]
{Independence of~$\boldsymbol{\ell}$ of Monodromy Groups}

\author{CheeWhye Chin}

\address{Department of Mathematics\\
	 University of California\\
	 Berkeley, CA 94720, U.S.A.}

\curraddr{The Broad Institute -- MIT\\
          320 Charles Street\\
	  Cambridge, MA 02141, U.S.A.}

\email{cheewhye@math.berkeley.edu\\
       cheewhye@mit.edu }


\dedicatory{Dedicated to Nicholas~M.~Katz
            on his 60th birthday}

\date{March 25, 2004}




%
%

\begin{abstract}

Let $X$ be a smooth curve
over a finite field of characteristic $p$,
let $E$ be a number field,
and let $\Lbold = \{\Lcal_\lambda\}$ be
an $E$-compatible system of lisse sheaves
on the curve $X$.
For each place $\lambda$ of $E$
not lying over~$p$,
the $\lambda$-component of the system $\Lbold$
is a lisse $E_\lambda$-sheaf $\Lcal_\lambda$ on $X$,
whose associated arithmetic monodromy group
is an algebraic group
over the local field $E_\lambda$.
We use
Serre's theory of Frobenius tori
and
Lafforgue's proof of Deligne's conjecture
to show that
when the $E$-compatible system $\Lbold$ is
semisimple and pure of some integer weight,
the isomorphism type of
the identity component of
these monodromy groups
is ``independent of~$\lambda$''.
More precisely,
after replacing $E$ by a finite extension,
there exists
a connected split reductive algebraic group
$G_0$ over the number field $E$
such that
for every place $\lambda$ of $E$
not lying over~$p$,
the identity component of
the arithmetic monodromy group of $\Lcal_\lambda$
is isomorphic to
the group $G_0$ with coefficients extended
to the local field $E_\lambda$.
\end{abstract}

\maketitle



\section{Introduction}
\label{sect:intro}

Let $X$ be a smooth curve
over a finite field of characteristic~$p$.
Let $E$ be a number field,
and consider
an $E$-compatible system $\Lbold = \{\Lcal_\lambda\}$
of lisse sheaves on $X$.
This means that
for every place $\lambda$ of $E$
not lying over~$p$,
we are given
a lisse $E_\lambda$-sheaf $\Lcal_\lambda$ on $X$,
and these lisse sheaves
are $E$-compatible with one another
in the sense that,
for every closed point $x$ of $X$,
the polynomial
$$ \det(1-T\,\Frob_x, \Lcal_\lambda)
$$
has coefficients in $E$
and is independent of the place $\lambda$.

Let $\etabar \rightarrow X$ be a geometric point of $X$.
For every place $\lambda$ of $E$
not lying over~$p$,
the lisse $E_\lambda$-sheaf $\Lcal_\lambda$,
which forms the $\lambda$-component of
the system $\Lbold$,
has its corresponding
monodromy $E_\lambda$-representation
$$ [\Lcal_\lambda] : \pi_1(X,\etabar) \rightarrow
   \GL({\Lcal_\lambda}_\etabar),
$$
and its corresponding
{\sl arithmetic monodromy group\/}
$\Garith(\Lcal_\lambda,\etabar)$,
defined as the Zariski closure of
the image of $\pi_1(X,\etabar)$
under $[\Lcal_\lambda]$.
The group $\Garith(\Lcal_\lambda,\etabar)$
is an algebraic group
over the local $\ell$-adic field $E_\lambda$,
and it is given with
a faithful tautological
$E_\lambda$-rational representation
$$ \sigma_\lambda:
   \Garith(\Lcal_\lambda,\etabar)
   \lhook\joinrel\rightarrow
   \GL({\Lcal_\lambda}_\etabar).
$$
The philosophy of motives in algebraic geometry
leads one to suspect that
if our $E$-compatible system $\Lbold = \{\Lcal_\lambda\}$
is ``motivic'',
then the collection
$\{\Garith(\Lcal_\lambda,\etabar), \sigma_\lambda\}$
of $\ell$-adic monodromy groups
and their tautological representations
should be
``independent of~$\ell$''
in a suitable sense.
This philosophy can be made more precise
in the form of the following conjecture,
which we are going to address in this paper.

\begin{conjecture}
\label{conj:main conj}
Let $X$ be a smooth curve
over a finite field of characteristic~$p$.
Let $E$ be a number field,
and let $\Lbold = \{\Lcal_\lambda\}$ be
an $E$-compatible system of lisse sheaves
on the curve $X$.
Assume that
the $E$-compatible system $\Lbold$ is
semisimple and
pure of weight~$w$ for some integer~$w$.
\begin{itemize}
\item[(i)]
There exist
a finite extension $F$ of $E$
and an algebraic group $G$ over the number field $F$
such that
for every place $\lambda\in|F|_{\ne p}$ of $F$
not lying over~$p$,
writing $\lambda$ also for its restriction to $E$,
the $F_\lambda$-groups
$$ G \otimes_F F_\lambda
   \qquad\text{and}\qquad
   \Garith(\Lcal_\lambda,\etabar) \otimes_{E_\lambda} F_\lambda
\qquad\text{are isomorphic}.
$$
\item[(ii)]
Assume that~(i) holds.
After replacing $F$ by a further finite extension,
there exists
an $F$-rational representation $\sigma$ of $G$
such that
for every place $\lambda\in|F|_{\ne p}$ of $F$
not lying over~$p$,
writing $\lambda$ also for its restriction to $E$,
and identifying
the $F_\lambda$-groups
$G \otimes_F F_\lambda$
and
$\Garith(\Lcal_\lambda,\etabar) \otimes_{E_\lambda} F_\lambda$
via an isomorphism given in~(i),
the $F_{\lambda}$-rational representations
$$ \sigma \otimes_F F_\lambda
   \qquad\text{and}\qquad
   \sigma_\lambda \otimes_{E_\lambda} F_\lambda
\qquad\text{are isomorphic}.
$$
\end{itemize}
\end{conjecture}

In the situation of~(i),
one could be more optimistic
and ask for
an algebraic group $G$ over the number field $E$ itself,
i.e.~without first passing to a finite extension $F$.
I do not know
if that would constitute ``asking for too much'',
because the philosophy of motives
does not predict that
the group $G$ can be defined
over the given field $E$.
However,
in the situation of~(ii),
even if the group $G$
is defined over $E$,
one would in general
need to pass to a finite extension
in order to find the representation $\sigma$.

A conjecture similar to~(\ref{conj:main conj})
for compatible systems of cohomology sheaves
has been studied previously
by M.~Larsen and R.~Pink;
see \cite{LarsenPink-Abvar}, Conj.~5.1.
They have also addressed
the related situation
in the ``abstract'' setting;
see \cite{LarsenPink-Indepl}, especially~\S9.

In the situation of~(\ref{conj:main conj}),
for each place $\lambda\in |E|_{\ne p}$ of $E$
not lying over $p$,
the group $\Garith(\Lcal_\lambda,\etabar)$
sits in a short exact sequence
$$ 1 \rightarrow \Garith(\Lcal_\lambda,\etabar)^0
     \rightarrow \Garith(\Lcal_\lambda,\etabar)
     \rightarrow \Garith(\Lcal_\lambda,\etabar) \ / \ \Garith(\Lcal_\lambda,\etabar)^0
     \rightarrow 1.
$$
The restriction of $\sigma_\lambda$
to $\Garith(\Lcal_\lambda,\etabar)^0$
will be denoted by the same symbol.
By a result of J.-P.~Serre,
one knows that
the finite groups
$\Garith(\Lcal_\lambda,\etabar) \ / \ \Garith(\Lcal_\lambda,\etabar)^0$
are independent of~$\lambda$.
More precisely:

\begin{theorem}[J.-P.~Serre]
\label{thm:Serre indep kernel arith}
{\rm (See~\cite{Serre-LettresRibet}, Th.~on p.~15.)}
With the notation and hypotheses of~(\ref{conj:main conj}),
the kernel of
the surjective homomorphism
$$ \pi_1(X,\etabar) \xrightarrow{[\Lcal_\lambda]}
   \Garith(\Lcal_\lambda,\etabar) \xrightarrow\joinrel\rarrow
   \Garith(\Lcal_\lambda,\etabar) \ / \ \Garith(\Lcal_\lambda,\etabar)^0
$$
is the same
open subgroup of $\pi_1(X,\etabar)$
for every place $\lambda\in|E|_{\ne p}$ of $E$
not lying over~$p$.
\end{theorem}

There are two proofs of this
in \cite{Serre-LettresRibet}, pp.~15--20;
a third proof can be found
in \cite{LarsenPink-Indepl}, Prop.~6.14.
As a consequence:

\begin{corollary}
With the notation and hypotheses of~(\ref{conj:main conj}),
there exists a finite group $\Gamma$
such that
for every place $\lambda\in|E|_{\ne p}$ of $E$
not lying over~$p$,
the finite groups
$$ \Gamma	\qquad\text{and}\qquad
   \Garith(\Lcal_\lambda,\etabar) \ / \ \Garith(\Lcal_\lambda,\etabar)^0
\qquad\text{are isomorphic}.
$$
\end{corollary}

A simplified version of our main result
asserts
the ``independence of~$\lambda$''
of the identity component
$\Garith(\Lcal_\lambda,\etabar)^0$
of the monodromy groups $\Garith(\Lcal_\lambda,\etabar)$,
and of their
tautological representations $\sigma_\lambda$:

\begin{theorem}
\label{thm:main arith}
Assume the notation and hypotheses
of~(\ref{conj:main conj}).
\begin{itemize}
\item[(i)]
There exist
a finite extension $F$ of $E$
and a connected split reductive algebraic group $G_0$
over the number field $F$
such that
for every place $\lambda\in|F|_{\ne p}$ of $F$
not lying over~$p$,
writing $\lambda$ also for
its restriction to $E$,
the connected $F_\lambda$-algebraic groups
$$ G_0 \otimes_F F_\lambda
   \qquad\text{and}\qquad
   \Garith(\Lcal_\lambda,\etabar)^0 \otimes_{E_\lambda} F_\lambda
\qquad\text{are isomorphic}.
$$
\item[(ii)]
There exists
an $F$-rational representation $\sigma_0$ of $G_0$
such that
for every place \linebreak
 $\lambda\in|F|_{\ne p}$ of $F$
not lying over~$p$,
writing $\lambda$ also for its restriction to $E$,
and identifying
the $F_\lambda$-groups
$G_0 \otimes_F F_\lambda$
and
$\Garith(\Lcal_\lambda,\etabar)^0 \otimes_{E_\lambda} F_\lambda$
via an isomorphism given in~(i),
the $F_{\lambda}$-rational representations
$$ \sigma_0 \otimes_F F_\lambda
   \qquad\text{and}\qquad
   \sigma_\lambda \otimes_{E_\lambda} F_\lambda
\qquad\text{are isomorphic}.
$$
\end{itemize}
\end{theorem}

Note that
the isomorphism in~(i)
between the two groups
is not unique,
nor is there a canonical one,
but the hypotheses of the theorem
do allow us
to rigidify the situation to some extent;
we refer to~(\ref{cor:stronger main thm arith})
for the precise statement.

Suppose the curve $X$ is 
geometrically connected over the base field $k$.
Let $\kbar$ be
the algebraic closure of $k$ in $\kappa(\etabar)$,
and regard $\etabar$ also as
a geometric point of $X\otimes_k\kbar$.
Then for every place $\lambda$ of $E$
not lying over~$p$,
we may also consider
the {\sl geometric monodromy group\/}
$\Ggeom(\Lcal_\lambda,\etabar)$
of the lisse $E_\lambda$-sheaf $\Lcal_\lambda$,
i.e.~the Zariski closure of
the image of $\pi_1(X\otimes_k\kbar,\etabar)$
under $[\Lcal_\lambda]$;
we write $\sigma'_\lambda$ for
its faithful tautological
$E_\lambda$-rational representation.
Thanks to Larsen and Pink,
one has the analogue of~(\ref{thm:Serre indep kernel arith}):

\begin{theorem}[M.~Larsen and R.~Pink]
\label{thm:LarsenPink indep kernel geom}
{\rm (See~\cite{LarsenPink-Abvar}, Prop.~2.2.)}
With the notation and hypotheses of~(\ref{conj:main conj}),
the kernel of
the surjective homomorphism
$$ \pi_1(X\otimes_k\kbar,\etabar) \xrightarrow{[\Lcal_\lambda]}
   \Ggeom(\Lcal_\lambda,\etabar) \xrightarrow\joinrel\rarrow
   \Ggeom(\Lcal_\lambda,\etabar) \ / \ \Ggeom(\Lcal_\lambda,\etabar)^0
$$
is the same
open subgroup of $\pi_1(X\otimes_k\kbar,\etabar)$
for every place $\lambda\in|E|_{\ne p}$ of $E$
not lying over~$p$.
\end{theorem}

From a result of P.~Deligne
(cf.~\cite{Deligne-WeilII}, Cor.~1.3.9),
one infers that
the identity component
$\Ggeom(\Lcal_\lambda,\etabar)^0$
of the geometric mono\-dromy group
is equal to
the derived subgroup
$\Der(\Garith(\Lcal_\lambda,\etabar)^0)$
of the identity component
of the arithmetic monodromy group.
Consequently,
our main theorem~(\ref{thm:main arith})
implies the following:

\begin{theorem}
\label{thm:main geom}
Assume the notation and hypotheses
of~(\ref{conj:main conj}).
\begin{itemize}
\item[(i)]
There exist
a finite extension $F$ of $E$
and a connected split semisimple algebraic group $G'_0$
over the number field $F$
such that
for every place $\lambda\in|F|_{\ne p}$ of $F$
not lying over~$p$,
writing $\lambda$ also for
its restriction to $E$,
the connected $F_\lambda$-algebraic groups
$$ G'_0 \otimes_F F_\lambda
   \qquad\text{and}\qquad
   \Ggeom(\Lcal_\lambda,\etabar)^0 \otimes_{E_\lambda} F_\lambda
\qquad\text{are isomorphic}.
$$
More precisely,
one can take $G'_0$
to be the derived subgroup of
the group $G_0$
in~(\ref{thm:main arith})(i).
\item[(ii)]
There exists
an $F$-rational representation $\sigma'_0$ of $G'_0$
such that
for every place \linebreak
 $\lambda\in|F|_{\ne p}$ of $F$
not lying over~$p$,
writing $\lambda$ also for its restriction to $E$,
and identifying
the $F_\lambda$-groups
$G'_0 \otimes_F F_\lambda$
and
$\Ggeom(\Lcal_\lambda,\etabar)^0 \otimes_{E_\lambda} F_\lambda$
via an isomorphism given in~(i),
the $F_{\lambda}$-rational representations
$$ \sigma'_0 \otimes_F F_\lambda
   \qquad\text{and}\qquad
   \sigma'_\lambda \otimes_{E_\lambda} F_\lambda
\qquad\text{are isomorphic}.
$$
More precisely,
one can take $\sigma'_0$
to be the restriction to $G'_0$
of the representation $\sigma_0$
in~(\ref{thm:main arith})(ii).
\end{itemize}
\end{theorem}

A weaker form of~(\ref{thm:main geom}),
obtained earlier by Larsen and Pink
(see~\cite{LarsenPink-Abvar}, Th.~2.4),
asserts that
after scalar extensions
from the various $\ell$-adic local fields
to a common algebraically closed field
(such as $\bC$),
the various $\Ggeom(\Lcal_\lambda,\etabar)^0$
become isomorphic to one another,
and that the same goes for 
the representations $\sigma'_\lambda$
under the additional hypothesis that
the compatible system is
geometrically absolutely irreducible.
Their results were established
using their theorem
(cf.~\cite{LarsenPink-DetermineRep}, Th.~1 and Th.~2)
to the effect that
a connected semisimple algebraic group
over an algebraically closed field
of characteristic~$0$
is determined
up to (non-canonical) isomorphism
by the {\it dimension data\/} of the group.
As we shall explain in a moment,
we adopt a completely different approach
in this paper
by making full use of
the motivic origin
of compatible systems on a curve.

The results of~(\ref{thm:main arith})
and~(\ref{thm:main geom})
(and for that matter,
 those of~(\ref{thm:Serre indep kernel arith})
 and~(\ref{thm:LarsenPink indep kernel geom})
 as well)
also hold
when $X$ is any
irreducible normal variety
of finite type over a finite field
of charac\-teristic~$p$.
This is seen by
reducing to the case of curves
---
for instance,
by using space-filling curves
(cf.~\cite{Katz-SpaceFillCurves}, Th.~8 and Lemma~6).

The proof of~(\ref{thm:main arith})
is given in \S\ref{sect:Indep of ell};
we outline the strategy here
as a guide to the organization of this paper.
For each $\lambda\in|E|_{\ne p}$,
write $G_{E_\lambda}$ for
$\Garith(\Lcal_\lambda,\etabar)$.
To fix ideas,
let us assume (only for simplicity) that 
$G_{E_\lambda}$ is connected
for {\it one\/} $\lambda\in|E|_{\ne p}$;
by~(\ref{thm:Serre indep kernel arith}),
the same is then true
for {\it every\/} $\lambda\in|E|_{\ne p}$.
We want to show that,
allowing $E$ to be replaced
by a finite extension,
the $G_{E_\lambda}$ for various $\lambda$
are all obtained by scalar extensions
from a common group $G_0$ over $E$.

The key result
which lies at the core of
our whole argument
is an extension~(\ref{thm:indeplalg})
of the fundamental theorem~(\ref{thm:Lafforgue 3})
of L.~Lafforgue;
it allows us to exploit
the hypothesis of compatibility
(cf.~(\ref{prop:mu to lambda})
and~(\ref{thm:bijection Irr}))
to establish that
the Grothendieck rings
of the various $G_{E_\lambda}$
are all isomorphic
in a way
which identifies the irreducible representations
and which respects
the character
of each irreducible representation.
Roughly speaking,
this means that
the connected reductive groups $G_{E_\lambda}$
all ``have the same representation theory''.
Making use of the tight connection
(cf.~(\ref{thm:determine gp}))
between the representation theory
and the structure theory
of connected reductive groups,
we can then conclude that
the various $G_{E_\lambda}$
do come from a common source $G_0$,
{\it provided that\/}
they contain maximal tori
coming from
a common {\it split\/} torus $T_0$
over $E$.
Granting this,
the fact that
the representations $\sigma_\lambda$ of $G_{E_\lambda}$
are all obtained by scalar extensions
from a common representation $\sigma_0$ of $G_0$
is then an easy consequence
(cf.~(\ref{proof:main(ii)}))
of our constructions.

To get the torus $T_0$,
we appeal to Serre's theory of Frobenius tori
(cf.~(\ref{thm:Serre max torus}),
     (\ref{cor:max dim torus}));
accordingly,
we can find a Frobenius element
which generates an $E$-torus
whose $\ell$-adic scalar extensions
yield maximal tori
in the various $G_{E_\lambda}$;
hence,
if we make an extension of $E$
to split this torus,
we obtain the $T_0$ we want.
However, in order to apply Serre's results,
we have to verify that
the lisse sheaves $\Lcal_\lambda$
satisfy certain hypotheses.
All but one of these hypotheses
are provided by
another theorem~(\ref{thm:Lafforgue 12abc})
of Lafforgue;
the ``missing'' hypothesis
(concerning boundedness in denominator)
is a result~(\ref{thm:2d})
which we establish in \S\ref{sect:boundedness denom}.

In the course of this work,
I have benefited tremendously
from discussions with
Pierre Deligne, Johan de~Jong,
Nicholas Katz and Laurent Lafforgue.
It is a pleasure
to acknowledge
my intellectual debts
to all of them.
I am also grateful
to the anonymous referees,
whose valuable suggestions have helped
improve and clarify the exposition here.


\section{Definitions and notation}
\label{sect:Def and Notation}

\subsection{}
Let $E$ be a number field.
If $p$ is a prime number,
let $|E|_p$ denote
the finite set of valuations of $E$ lying over~$p$,
and let
$$ |E|_{\ne p} := \bigcup_{\text{$\ell$ prime $\ne p$}} |E|_\ell
$$
denote the set of valuations of $E$
not lying over~$p$.
Let $|E|_\infty$ denote
the finite set of archimedean absolute values of $E$.
We regard each $\nu\in |E|_p$
as a homomorphism
$$ \nu : E^\times \rightarrow \bQ,
	\qquad\text{normalized so that
	$\nu(p) = 1$},
$$
and we regard each $|\cdot|\in |E|_\infty$
as a homomorphism
$$ |\cdot| : E^\times \rightarrow \bR^\times_{> 0},
	\qquad\text{normalized so that
	$|n| = n$ for every $n\in\bZ_{>0}$}.
$$
With these normalizations,
the valuations and absolute values
are compatible with
passing to finite extensions of $E$.

Let $\alpha\in \bQbar^\times$ be
a nonzero algebraic number.
Let $p$ be a prime number,
and let $q$ be a positive power of $p$.
\begin{enumerate}
\item
Let $w\in\bZ$ be an integer;
we say that
$\alpha$ is {\sl pure of weight~$w$ with respect to~$q$\/}
if for every 
archimedean absolute value
$|\cdot| \in |\bQ(\alpha)|_\infty$ of $\bQ(\alpha)$,
one has
$$ |\alpha| = q^{w/2}.
$$
\item
We say that
$\alpha$ is {\sl plain of characteristic~$p$\/}
if $\alpha$ is an $\ell$-adic unit
for every prime $\ell\ne p$;
in other words,
for every
non-archimedean valuation
$\lambda \in |\bQ(\alpha)|_{\ne p}$ of $\bQ(\alpha)$
not lying over~$p$,
one has
$$ \lambda(\alpha) = 0.
$$
\item
Let $C \ge 0$ be a real number;
we say that
$\alpha$ is {\sl $C$-bounded in valuation with respect to~$q$\/}
if for every 
non-archimedean valuation
$\nu \in |\bQ(\alpha)|_p$ of $\bQ(\alpha)$
lying over~$p$,
one has
$$ \left| \dfrac{\nu(\alpha)}{\nu(q)} \right| \le C.
$$
\item
Let $D > 0$ be an integer;
we say that
$\alpha$ is {\sl $D$-bounded in denominator with respect to~$q$\/}
if for every 
non-archimedean valuation
$\nu \in |\bQ(\alpha)|_p$ of $\bQ(\alpha)$
lying over~$p$,
one has
$$ \dfrac{\nu(\alpha)}{\nu(q)} \in \tfrac{1}{D}\,\bZ.
$$
\end{enumerate}

\subsection{}
Let $\ell$ be a prime number,
and let $\Lambda$ be
an $\ell$-adic field
(i.e.~$\Lambda$ is
an algebraic extension of $\bQell$).
Let $X$ be
a connected normal scheme
of finite type over $\Spec(\bZ[1/\ell])$,
and let $\etabar \rightarrow X$ be
a geometric point of $X$.
If $\Lcal$ is a lisse $\Lambda$-sheaf on $X$,
we let
$$ [\Lcal]: \pi_1(X,\etabar) \rightarrow \GL(\Lcal_\etabar)
$$
denote the corresponding
continuous monodromy $\Lambda$-representation
of the arithmetic fundamental group
$\pi_1(X,\etabar)$ of $X$.
The {\sl arithmetic monodromy group\/}
$\Garith(\Lcal,\etabar)$
of $\Lcal$
is the Zariski closure of
the image of $\pi_1(X,\etabar)$ in $\GL(\Lcal_\etabar)$
under the monodromy representation $[\Lcal]$;
it is a linear algebraic group over $\Lambda$.

\subsection{}
Let $|X|$ denote
the set of closed points of $X$.
For each $x\in|X|$,
choose an algebraic geometric point
$\xbar \rightarrow x \in X$ of $X$
localized at $x$.
The absolute Galois group $\Gal(\kxbar/\kx)$
is a free profinite group
generated by
the {\sl geometric Frobenius at $x$\/},
$\Frob_x := \Frob_\kx \in \Gal(\kxbar/\kx)$.
If $\Lcal$ is
a lisse $\Lambda$-sheaf on $X$,
we denote its restriction to $x$
by $\Lcal(x)$,
and we write
$$ [\Lcal(x)] : \Gal(\kxbar/\kx) \rightarrow \GL(\Lcal_\xbar)
$$
for the monodromy $\Lambda$-representation
of $\Lcal(x)$.
Thus, $[\Lcal(x)](\Frob_x)$
is an element of $\GL(\Lcal_\xbar)$.

Let $\Lambdabar$ be
an algebraically closed extension of $\Lambda$.
A lisse $\Lambda$-sheaf $\Lcal$ on $X$ is
\begin{itemize}
\item[--]
{\sl algebraic\/},
\item[--]
resp.\ 
{\sl pure of weight~$w$\/}
(for an integer~$w$),
\item[--]
resp.\ 
{\sl plain of characteristic~$p$\/}
(for a prime number~$p$),
\item[--]
resp.\ 
{\sl $C$-bounded in valuation\/}
(for a real number~$C\ge 0$),
\item[--]
resp.\ 
{\sl $D$-bounded in denominator\/}
(for an integer~$D > 0$)
\end{itemize}
if for every closed point $x\in|X|$,
and for every eigenvalue $\alpha\in\Lambdabar^\times$
of $[\Lcal(x)](\Frob_x)$,
\begin{itemize}
\item[--]
$\alpha$ is an algebraic number,
\item[--]
resp.\ 
$\alpha$ is pure of weight~$w$
with respect to $\#\kx$,
\item[--]
resp.\ 
$\alpha$ is plain of characteristic~$p$,
\item[--]
resp.\ 
$\alpha$ is $C$-bounded in valuation
with respect to $\#\kx$,
\item[--]
resp.\ 
$\alpha$ is $D$-bounded in denominator
with respect to $\#\kx$.
\end{itemize}
These properties of lisse $\Lambda$-sheaves
are stable under
passage to subquotients
and formation of extensions.
Let $a: X' \rightarrow X$ be a morphism of finite type;
if $\Lcal$ on $X$ has one of these properties,
then so does $a^*\Lcal$ on $X'$,
and the converse holds when $a$ is surjective.

\subsection{}
Let $X$ be a scheme of finite type
over a finite field of characteristic~$p$.
Let $E$ be a number field.
An {\sl $E$-compatible system $\Lbold$ on $X$\/}
is a collection $\Lbold = \{\Lcal_\lambda\}$
of lisse sheaves on $X$
indexed by the set $|E|_{\ne p}$
of places of $E$ not lying over~$p$,
where, for every such place
$\lambda\in |E|_{\ne p}$ of $E$,
$\Lcal_\lambda$ is
a lisse $E_\lambda$-sheaf on $X$,
and this collection of lisse sheaves
satisfies the condition that
they are $E$-compatible with one another:
i.e.~for every closed point $x\in|X|$ of $X$,
and for every place $\lambda\in|E|_{\ne p}$ of $E$
not lying over~$p$,
the polynomial
$$ \det(1-T\,\Frob_x, \Lcal_\lambda)
$$
has coefficients in $E$
and is independent of the place $\lambda$.

We say that
the $E$-compatible system $\Lbold$ is
{\sl pure of weight~$w$\/}
(for an integer~$w$)
if for any/every $\lambda\in|E|_{\ne p}$,
the lisse $E_\lambda$-sheaf $\Lcal_\lambda$
is pure of weight~$w$.
We say that $\Lbold$ is
{\sl absolutely irreducible\/}, 
resp.~{\sl semi\-simple\/},
if for each $\lambda\in|E|_{\ne p}$,
the lisse $E_\lambda$-sheaf $\Lcal_\lambda$
has the corresponding property
(i.e.~is irreducible over an algebraic closure of $E_\lambda$,
resp.~is semisimple over $E_\lambda$).

Note that
if $\Lbold = \{\Lcal_\lambda\}$ is
an $E$-compatible system,
then for every $\lambda\in|E|_{\ne p}$,
the lisse $E_\lambda$-sheaf $\Lcal_\lambda$
is necessarily
plain of characteristic~$p$.
We will see later (cf.~(\ref{thm:indeplalg})) that
this necessary condition is also sufficient
for a lisse $\ell$-adic sheaf
to extend to
an $E$-compatible system
for some number field $E$.


\section{Boundedness of denominator}
\label{sect:boundedness denom}

Let us first recall
the following fundamental result,
first conjectured by P.~Deligne,
and now established by L.~Lafforgue.

\begin{theorem}[L.~Lafforgue]
\label{thm:Lafforgue 12abc}
{\rm (See~\cite{Lafforgue-LanglandsCorresp}, Th\'eor\`eme~VII.6~(i--iv).)}
Let $X$ be a smooth curve
over a finite field of characteristic~$p$.
Let $\ell \ne p$ be a prime number,
and let $\Lcal$ be a lisse $\bQellbar$-sheaf on $X$
which is irreducible,
of rank~$r$,
and whose determinant is of finite order.
Then
there exists a number field $E \subset \bQellbar$
such that
the lisse sheaf $\Lcal$ on $X$ is
\begin{enumerate}
\item[(1)]
$E$-rational (hence algebraic);
\item[(2a)]
pure of weight~$0$;
\item[(2b)]
plain of characteristic~$p$; and
\item[(2c)]
$C$-bounded in valuation,
where $C$ may be taken to be $(r-1)^2/r$.
\end{enumerate}
\end{theorem}

In this section, we prove the following assertion~(2d),
which complements~(\ref{thm:Lafforgue 12abc}) above.

\begin{theorem}
\label{thm:2d}
With the notation and hypotheses of~(\ref{thm:Lafforgue 12abc}),
there exists an integer~$D > 0$
such that
the lisse sheaf $\Lcal$ on $X$ is
\begin{enumerate}
\item[(2d)]
$D$-bounded in denominator.
\end{enumerate}
\end{theorem}

\begin{proof}[Proof of~(\ref{thm:2d}), prelude]
We begin our argument
as Lafforgue does
for his proof of assertion~(2c)
in~(\ref{thm:Lafforgue 12abc});
see \cite{Lafforgue-LanglandsCorresp}, Th.~VII.6, D\'emo.,
pp.~198--200, especially parts~(v) and~(iv).
Accordingly,
we consider
the lisse sheaf
$$ {q_1}^*\Lcal\otimes {q_2}^*\Lcal^{\vee}(1-r)
	\qquad\text{on the surface $X \times X$},
$$
where $q_1$ and $q_2$ denote
the two projections from $X\times X$ onto $X$.
The fundamental fact we need to know
about this lisse sheaf
is that
it ``appears'' in
the $\ell$-adic cohomology of a certain stack
over the surface.
More precisely,
there exists
a stack $\Xfrak$ over $X\times X$
---
namely, $\Chtbar/a^\bZ$
in the notation of \cite{Lafforgue-LanglandsCorresp}
---
such that if
$$ f: \Xfrak \rightarrow X\times X
$$
denotes the structural morphism,
the semisimple\footnote{\ 
Cf.~\cite{Lafforgue-LanglandsCorresp},
remark after D\'ef.~VI.14.}
lisse sheaf
${q_1}^*\Lcal\otimes {q_2}^*\Lcal^{\vee}(1-r)$
on $X \times X$
occurs as a direct summand
in the semisimplified
lisse\footnote{\ 
Cf.~\cite{Lafforgue-LanglandsCorresp},
Chap.~VI,~\S2a, especially p.~165, third paragraph.}
cohomology sheaf
$(\tR^{2r-2} f_! (\bQellbar))^\semisimp$
of $\Xfrak$
over the surface
$X\times X$;
for the justification of this key assertion,
see~\cite{Lafforgue-LanglandsCorresp}, Chap.~VI,~\S3,
especially the statements of
Lemme~VI.26 and Th.~VI.27.
For our present purpose
of proving~(\ref{thm:2d}),
we need not be concerned with
the precise definition and moduli interpretation
of the stack $\Xfrak$;
what will be important for us
is the following assertion
obtained from
\cite{Lafforgue-LanglandsCorresp}, Th.~V.2 and Lemme~A.3:
there is a finite extension $\bF_{q'}$ of $\bF_q$,
and a scheme $Z$
---
namely, $(\Chtbar/a^\bZ)^\gr\otimes_{\bF_q} \bF_{q'}$
in the notation of \cite{Lafforgue-LanglandsCorresp}
---
which is proper over the surface
$V := (X \times X)\otimes_{\bF_q} \bF_{q'}$,
such that
if $p: Z \rightarrow V$ denotes
the structural morphism,
the cohomology sheaf
$\tR^{2r-2} p_! (\bQellbar)$
is equal to
the pull-back of
$\tR^{2r-2} f_! (\bQellbar)$
to $V$.
The point to note is that
if $\Gcal$ denotes
the pull-back of
the lisse sheaf
${q_1}^*\Lcal\otimes {q_2}^*\Lcal^{\vee}(1-r)$
to $V$,
then every irreducible constituent of $\Gcal$
is also an irreducible constituent
of the lisse cohomology sheaf
$\tR^{2r-2} p_! (\bQellbar)$
of a proper {\it scheme\/} $Z$ over $V$.
We use this observation
to deduce the following key lemma.

\renewcommand{\qedsymbol}{}
\end{proof}

\begin{lemma}
\label{lemma:denom}
Let $\Fcal$ be
an irreducible constituent of
the lisse sheaf $\Gcal$ on $V$.
\begin{enumerate}
\item
\label{lemma:denom(a)}
There exist
\begin{enumerate}
\item
a finite universal homeomorphism
$\Vtilde \rightarrow V$,
\item
an open dense subscheme $\Utilde$ of $\Vtilde$,
and
\item
a scheme $Z_\Fcal$ over $\Utilde$
whose structural morphism
$\ptilde_\Fcal: Z_\Fcal \rightarrow \Utilde$
is projective and smooth,
\end{enumerate}
such that
on the scheme $\Utilde$,
the lisse sheaf $\Fcal|_\Utilde$
is an irreducible constituent
of the cohomology sheaf
$\tR^{2r-2} (\ptilde_\Fcal)_!(\bQellbar)$
(the latter is a lisse sheaf on $\Utilde$
since $\ptilde_\Fcal$ is proper and smooth).
\item
\label{lemma:denom(b)}
Let $U$ denote
the image of $\Utilde$ in $V$.
There exists
an integer~$D > 0$
such that
the algebraic lisse sheaf $\Fcal|_U$ on $U$
is $D$-bounded in denominator.
\end{enumerate}
\end{lemma}

\begin{proof}[Proof of (\ref{lemma:denom}, (\ref{lemma:denom(a)}))]
Let $\eta$ be
the generic point of $V$.
The irreducible lisse sheaf
$\Fcal$ on $V$
restricts to
an irreducible lisse sheaf
$\Fcal_\eta$ on $\eta$.
We have the proper morphism
$p_\eta: Z_\eta \rightarrow \eta$,
whose cohomology sheaf
$\tR^{2r-2} (p_\eta)_! (\bQellbar)$
admits $\Fcal_\eta$
as an irreducible subquotient.

We now construct
a proper hypercovering
$a: Z'_\bullet \rightarrow Z_\eta$
of $Z_\eta$
(cf.~\cite{Deligne-HodgeIII},~(6.2.5)),
but using
the alteration theorem of de~Jong
(cf.~\cite{deJong-Alterations}, Th.~4.1)
instead of resolution.
Accordingly,
each $Z'_i$
is a projective scheme over $\eta$,
whose structural morphism
$p_i: Z'_i \rightarrow \eta$
factors through some $\eta_i$
in such a way that
$p'_i: Z'_i \rightarrow \eta_i$
is smooth,
and $\kappa(\eta_i) \supseteq \kappa(\eta)$
is a finite and purely inseparable extension
(this factorization property
of the structural morphisms
can be derived
from \cite{deJong-Alterations}, Remark~4.3).
Since a proper hypercovering
is of cohomological descent,
the adjunction morphism
$\bQellbar \rightarrow \tR a_! \bQellbar$
is an isomorphism
in the ``derived'' category $\tD^b_c(Z_\eta,\bQellbar)$
of $\bQellbar$-sheaves on $Z_\eta$.
We can then see
from the resulting
spectral sequence for $\tR (p_\eta)_!$
that
$\Fcal_\eta$ is an irreducible subquotient
of the cohomology sheaf
$\tR^j (p_i)_! (\bQellbar)$
of some $Z'_i$,
for some cohomological degree~$j$.

Let $\Vtilde$ be
the normalization of $V$
in $\eta_i$ over $\eta$.
Since $V$ is itself normal
and $\eta_i \rightarrow \eta$ is
finite and purely inseparable,
the morphism $\Vtilde \rightarrow V$
is a finite universal homeomorphism.
The projective and smooth morphism
$p'_i: Z'_i \rightarrow \eta_i$
``spreads out'' to
a projective and smooth morphism
$\ptilde_\Fcal: Z_\Fcal \rightarrow \Utilde$
for some open dense subscheme $\Utilde$ of $\Vtilde$
(i.e.~the latter morphism exists,
having the former morphism
as its generic fiber).
From the previous paragraph,
it follows that 
$\Fcal|\Utilde$ is
an irreducible subquotient of
$\tR^j (\ptilde_\Fcal)_! (\bQellbar)$.
Using Deligne's theory of weights
(cf.~\cite{Deligne-WeilII})
to compare the weights
of the lisse sheaves involved,
we see that
the cohomological degree $j$ above
is in fact equal to $2r-2$.
\end{proof}

\begin{proof}[Proof of (\ref{lemma:denom}, \ref{lemma:denom(b)}))]
From part~(\ref{lemma:denom(a)}),
we have the morphism
$\ptilde_\Fcal : Z_\Fcal \rightarrow \Utilde$,
which is projective and smooth.
Let $d$ be the rank of the lisse sheaf
$\Hcal := \tR^{2r-2} (\ptilde_\Fcal)_! (\bQellbar)$
on $\Utilde$,
and let $D := d!$.
By part~(\ref{lemma:denom(a)}),
$\Hcal$ admits $\Fcal|_\Utilde$
as an irreducible subquotient.
We will show that
the lisse $\bQellbar$-sheaf $\Hcal$ on $\Utilde$
is $\bQ$-rational
and is $D$-bounded in denominator;
this will imply that
$\Fcal|_\Utilde$ is also $D$-bounded in denominator,
and thus the conclusion of part~(\ref{lemma:denom(b)})
will follow.

Let $u\in|\Utilde|$ be a closed point,
$k = \kappa(u)$ be its residue field,
and $\kbar$ be an algebraic closure of $k$.
We have to show that
the $\bQellbar$-representation $[\Hcal(u)]$
of $\Gal(\kbar/k)$
is $\bQ$-rational
and is $d!$-bounded in denominator
with respect to~$\#k$.
Let $Y$ be the fiber of $Z_\Fcal$ over $u$;
thus $Y$ is a projective smooth scheme
over the finite field $k$.
By the proper base change theorem,
$[\Hcal(u)]$ is isomorphic to
the representation of $\Gal(\kbar/k)$
on the $\ell$-adic cohomology
$\tH^{2r-2}(Y\otimes_k \kbar, \bQellbar)$
of $Y\otimes_k \kbar$.
The fact that
$[\Hcal(u)]$ is $\bQ$-rational
then follows from
the Weil Conjectures
as proved by Deligne
(cf.~\cite{Deligne-WeilI} or \cite{Deligne-WeilII});
in fact, one knows that
the ``inverse'' characteristic polynomial
$\det(1-T\,[\Hcal(u)](\Frob_u))$
has coefficients in $\bZ$.
Note that
$d$ is the degree of this polynomial.

To see that
$[\Hcal(u)]$ is $d!$-bounded in denominator
with respect to~$\#k$,
we appeal to the crystalline cohomology
$H := \tH_\cris^{2r-2}(Y/W(k))$
of $Y$ with respect to
the ring $W(k)$ of Witt vectors over $k$;
for a concise summary of
the assertions we need,
see \cite{Berthelot-Slopes},~\S\S1--2.
The Frobenius endomorphism of $Y$ relative to $k$,
$$ \Fr_{Y/k} : Y \rightarrow Y,
	\qquad\text{defined by $y \mapsto y^{\#k}$},
$$
induces an endomorphism
$(\Fr_{Y/k})^* : H \rightarrow H$
which gives $H$
the structure of
an $(\Fr_{Y/k})^*$-crystal.
Let $K$ be the fraction field of $W(k)$,
and let $\nu_K$ be the discrete valuation of $K$,
which we regard as a homomorphism
$$ \nu_K : K^\times \rightarrow \bQ,
	\qquad\text{normalized so that
		    $\nu_K(\#k) = 1$}.
$$
One has the ``inverse'' characteristic polynomial
$$ \det(1-T\,(\Fr_{Y/k})^*, H\otimes_{W(k)} K) \in K[T]
$$
of the endomorphism of $H\otimes_{W(k)} K$
induced by $(\Fr_{Y/k})^*$.
By the compatibility theorem of Katz-Messing
(cf.~\cite{KatzMessing}, Th.~1 and Cor.~1),
one knows that
the coefficients of this polynomial
lie in $\bZ$,
and that one has
$$ \det(1-T\,(\Fr_{Y/k})^*, H\otimes_{W(k)} K)
 = \det(1-T\,[\Hcal(u)](\Frob_u))
	\qquad\text{(equality in $\bZ[T]$)},
$$
and hence
the $(\Fr_{Y/k})^*$-isocrystal $H\otimes_{W(k)} K$
has $K$-dimension $d$.
The desired assertion that
$[\Hcal(u)]$ is
$d!$-bounded in denominator with respect to~$\#k$
is therefore equivalent to
the assertion that
the Newton polygon of the polynomial
$\det(1-T\,(\Fr_{Y/k})^*, H\otimes_{W(k)} K)$
with respect to the valuation $\nu_K$
has all its slopes lying in $\tfrac{1}{d!}\,\bZ$.
By a theorem of Yu.~I.~Manin
(cf.~\cite{Berthelot-Slopes}, Th.~1.3(ii)),
the slopes of this Newton polygon
are the slopes of
the $(\Fr_{Y/k})^*$-isocrystal $H\otimes_{W(k)} K$,
and so by
the classification theorem of J.~Dieudonn\'e
(cf.~\cite{Berthelot-Slopes}, Th.~1.3(i)),
these slopes are
of the form $r/s$
where $r\in\bZ$,
$s \in \{1,\ldots,d\}$,
and $(r,s) = 1$;
in particular,
we see that
these slopes all lie in $\tfrac{1}{d!}\,\bZ$,
which is what we want.
\end{proof}

\begin{proof}[Proof of~(\ref{thm:2d}), coda]
Applying~(\ref{lemma:denom})
to every irreducible constituent
of the lisse sheaf $\Gcal$ on $V$
and consolidating the conclusions thus obtained,
we infer that
there exist
an open dense subscheme $U'$ of $V$
and an integer~$D_0 > 0$
such that
the algebraic lisse sheaf $\Gcal|_{U'}$ on $U'$
is $D_0$-bounded in denominator.
Let $U$ denote
the image of $U'$
under the finite etale map $V \rightarrow X \times X$.
Then we see that
the lisse sheaf
$({q_1}^*\Lcal\otimes {q_2}^*\Lcal^{\vee}(1-r))|_U$ on $U$
is also $D_0$-bounded in denominator.
The complement of $U$ in $X \times X$
is a finite union of
divisors and closed points;
hence
$$ S := \left\{\  x \in |X| \ :\  
	  {q_1}^{-1}(x) \cap U = \emptyset
	  \quad\text{in $X \times X$}  \ \right\}
	\qquad\text{is a finite set}.
$$

Let $D = r\,D_0$.
We claim that
the algebraic lisse sheaf $\Lcal|_{X-S}$ on $X-S$
is $D$-bounded in denominator.
To see this,
let $x\in |X|-S$ be a closed point,
let $\alpha\in\bQellbar^\times$ be
an eigenvalue of $\Frob_x$
acting on $\Lcal|_{X-S}$,
and let $\nu\in|\bQ(\alpha)|_p$ be
a non-archimedean valuation of $\bQ(\alpha)$
lying over~$p$.
By the definition of $S$,
we can find
a closed point $u$ of $U$
with $q_1(u) = x$.
Set $y := q_2(u) \in |X|$,
and let
$e = [\kappa(u):\kx]$,
$f = [\kappa(u):\kappa(y)]$.
If $\beta\in\bQellbar^\times$
is an eigenvalue of $\Frob_y$ acting on $\Lcal$,
it follows that
$\alpha^e\,\beta^{-f}\,(\#\kappa(u))^{2r-2} \in \bQellbar^\times$
is an eigenvalue of $\Frob_u$
acting on $({q_1}^*\Lcal\otimes {q_2}^*\Lcal^{\vee}(1-r))|_U$;
so we have
$$ \dfrac{\nu(\alpha^e\,\beta^{-f}\,(\#\kappa(u))^{2r-2})}{\nu(\#\kappa(u))}
	 \in \tfrac{1}{D_0}\,\bZ,
$$
and hence
$$ \dfrac{e\,\nu(\alpha) - f\,\nu(\beta)}{\nu(\#\kappa(u))}
	 \in \tfrac{1}{D_0}\,\bZ.
$$
By hypothesis,
the determinant of $\Lcal$
is of finite order,
which implies that
the product of the $r$~many eigenvalues $\beta$
of $\Frob_y$ acting on $\Lcal$
is a root of unity;
it follows that
$$ \dfrac{r e\,\nu(\alpha)}{\nu(\#\kappa(u))}
 = \dfrac{r\,\nu(\alpha)}{\nu(\#\kx)}
	 \in \tfrac{1}{D_0}\,\bZ,
$$
which proves our claim.

Since $S$ is a finite set,
we may now replace $D$ by a multiple
to ensure that
the algebraic lisse sheaf $\Lcal$
on the whole of $X$
is $D$-bounded in denominator.
This completes the proof of~(\ref{thm:2d}).
\end{proof}

\begin{remark}
\label{rmk:generalization for 2abc}
Lafforgue has also shown
(see \cite{Lafforgue-LanglandsCorresp}, Prop.~VII.7) that,
as predicted by Deligne's conjecture,
assertions~(2a), (2b) and~(2c) of~(\ref{thm:Lafforgue 12abc})
generalize to the case
when $X$ is an arbitrary
normal variety of finite type
over a finite field of charac\-teristic~$p$.
The proof is by
reduction to the case of curves,
using the fact that
the assertions in~(\ref{thm:Lafforgue 12abc})
are {\it uniform\/} for
any lisse sheaf of a given rank~$r$
on any curve.
Unfortunately,
our proof of~(\ref{thm:2d})
does not produce
an expression for the integer~$D>0$
in terms of
just the rank $r$ of the lisse sheaf,
and so
we do not get
the desired generalization of assertion~(2d).
\end{remark}


\section{Compatible systems of Lisse sheaves}
\label{sect:compatible systems}

According to Deligne's conjecture
(cf.~\cite{Deligne-WeilII}, Conj.~(1.2.10)(v)),
on any normal connected scheme $X$
over a finite field,
an irreducible lisse $\bQellbar$-sheaf $\Lcal$
whose determinant $\det\Lcal$ is of finite order
should ``extend'' to
an absolutely irreducible
$E$-compatible system $\{\Lcal_\lambda\}$ of lisse sheaves
for some number field $E$.
Thanks to Lafforgue's proof
of the Langlands Correspondence
for $\GL_r$ over function fields,
this conjecture is now known to hold
when $X$ is a curve:

\begin{theorem}
\label{thm:Lafforgue 3}
{\rm (cf.~\cite{Lafforgue-LanglandsCorresp},~Th\'eor\`eme~VII.6~(v)
 and \cite{Chin-IndeplLafforgue}) }
Let $X$ be a smooth curve
over a finite field of characteristic~$p$.
Let $\ell \ne p$ be a prime number,
and let $\Lcal$ be
an irreducible lisse $\bQellbar$-sheaf on $X$
whose determinant is of finite order.
Then
there exists a number field $E \subset \bQellbar$
and an absolutely irreducible $E$-compatible system
$\Lbold = \{\Lcal_\lambda\}$ on $X$
which is $E$-compatible
with the given lisse $\bQellbar$-sheaf $\Lcal$;
i.e.~for every place $\lambda\in|E|_{\ne p}$ of $E$
not lying over~$p$,
there exists an absolutely irreducible
lisse $E_\lambda$-sheaf $\Lcal_\lambda$ on $X$
which is $E$-compatible
with the given lisse $\bQellbar$-sheaf $\Lcal$.
\end{theorem}

The purpose of this section
is to give a mild generalization~(\ref{thm:indeplalg})
of the above theorem,
in which
one replaces the assumption
``determinant is of finite order''
by a weaker hypothesis of plainness.
We first recall (cf.~(\ref{def:Tate twist}))
the construction of
rank~$1$ twisting lisse sheaves,
and use these to give
equivalent characterizations (cf.~(\ref{prop:twist}))
of algebraicity, purity and plainness
for an irreducible lisse $\bQellbar$-sheaf.

\subsection{}
\label{def:Tate twist}
Let $\bF$ be
a finite field of characteristic~$p$,
let $\ell\ne p$ be a prime number,
and let $\Lambda$ be an $\ell$-adic field.
For any $\ell$-adic unit
$\alpha\in\Lambda^\times$,
one has a well-defined
continuous homomorphism
$$ \Gal(\bFbar/\bF) \rightarrow \GL_1(\Lambda) = \Lambda^\times,
   \qquad\text{mapping $\Frob_{\bF}$ to $\alpha$},
$$
which corresponds to
a lisse $\Lambda$-sheaf of rank~$1$ on
the scheme $\Spec(\bF)$;
we let the symbol $\alpha^{\deg_{\bF}}$ denote
the pullback of this lisse $\Lambda$-sheaf
to any scheme $X$ of finite type over $\bF$.
If $x\in|X|$ is a closed point of $X$,
then one has
$$ \det(1-T\,\Frob_x, \alpha^{\deg_{\bF}})
 = 1 - T\,\alpha^{\deg_{\bF}(\kappa(x))}
	\qquad\text{(equality in $\Lambda[T]$)}.
$$
If $\alpha = 1/(\#\bF)$,
the resulting lisse sheaf is
the {\sl Tate twist}.
An important special case:
if $E$ is a number field
and $\alpha\in E^\times$
is {\it plain of charac\-teristic~$p$},
then the rank~$1$ lisse $E_\lambda$-sheaf
$\alpha^{\deg_{\bF}}$
is defined
for every place $\lambda\in|E|_{\ne p}$ of $E$
not lying over~$p$,
and these lisse sheaves together
form an $E$-compatible system.

\begin{proposition}
\label{prop:twist}
Let $X$ be a smooth curve
over a finite field $\bF$ of characteristic~$p$.
Let $\ell \ne p$ be a prime number,
and let $\Lcal$ be a lisse $\bQellbar$-sheaf on $X$
which is irreducible, of rank~$r$.
Let $w\in\bZ$ be an integer.
The following are equivalent:
\begin{enumerate}
\item
\label{item:twist(a)}
$\Lcal$ is algebraic
(resp.~pure of weight~$w$,
resp.~plain of characteristic~$p$);
\item
\label{item:twist(b)}
for some closed point $x\in|X|$,
the $\bQellbar$-representation $[\Lcal(x)]$
of $\Gal(\kxbar/\kx)$ is algebraic
(resp.~pure of weight~$w$ with respect to $\#\kx$,
resp.~plain of characteristic~$p$);
\item
\label{item:twist(c)}
$\det\Lcal$ is algebraic
(resp.~pure of weight~$r\,w$,
resp.~plain of characteristic~$p$);
\item
\label{item:twist(d)}
for some closed point $x\in|X|$,
the $\bQellbar$-representation $[\det\Lcal(x)]$
of $\Gal(\kxbar/\kx)$ is algebraic
(resp.~pure of weight~$r\,w$,
resp.~plain of characteristic~$p$);
\item
\label{item:twist(e)}
there exist
an integer $n>0$
and an element $\beta\in\bQellbar^\times$
which is algebraic
(resp.~pure of weight~$n\,r\,w$ with respect to $q$,
resp.~plain of characteristic~$p$),
such that
$(\det\Lcal)^{\otimes n}$ is isomorphic to
$\beta^{\deg_{\bF}}$;
\item
\label{item:twist(f)}
there exists an element $\alpha\in\bQellbar^\times$
which is algebraic
(resp.~pure of weight~$w$ with respect to $q$,
resp.~plain of characteristic~$p$),
such that
$\Lcal\otimes \alpha^{\deg_{\bF}}$
is an irreducible lisse $\bQellbar$-sheaf
whose determinant is of finite order.
\end{enumerate}
\end{proposition}

\begin{proof}
It is clear that
$\text{(\ref{item:twist(a)})}\Rightarrow
 \text{(\ref{item:twist(b)})}\Rightarrow
 \text{(\ref{item:twist(d)})}$
and that
$\text{(\ref{item:twist(a)})}\Rightarrow
 \text{(\ref{item:twist(c)})}\Rightarrow
 \text{(\ref{item:twist(d)})}$.
By a result of Deligne
(cf.~\cite{Deligne-WeilII}, Prop.~(1.3.4)(i)),
one knows that
for some $n\in\bZ_{>0}$,
the $n$-th tensor power of $\det\Lcal$
is geometrically constant,
which means that
its monodromy representation
$[(\det\Lcal)^{\otimes n}]$
factors as
$$ \pi_1(X,\etabar)
 \rightarrow \Gal(\bFbar/\bF)
 \rightarrow \bQellbar^\times;
$$
hence $(\det\Lcal)^{\otimes n}$ is isomorphic to
$\beta^{\deg_{\bF}}$
where $\beta\in\bQellbar^\times$
is the image of
the geometric Frobenius
$\Frob_{\bF}\in\Gal(\bFbar/\bF)$;
thus
$\text{(\ref{item:twist(d)})}\Rightarrow
 \text{(\ref{item:twist(e)})}$.
Now let $\alpha := 1/{\root{nr}\of{\beta}} \in\bQbar^\times$
be an \linebreak
$nr$-th root of $1/\beta$ in $\bQbar$;
the determinant of
the lisse sheaf $\Lcal\otimes\alpha^{\deg_{\bF}}$
is then of finite order dividing~$n$,
whence
$\text{(\ref{item:twist(e)})}\Rightarrow
 \text{(\ref{item:twist(f)})}$.
Finally,
$\text{(\ref{item:twist(f)})}\Rightarrow
 \text{(\ref{item:twist(a)})}$
follows directly
from parts~(1), (2a) and~(2b) of
Lafforgue's theorem~(\ref{thm:Lafforgue 12abc}).
\end{proof}

\begin{remark}
The equivalent characterizations in~(\ref{prop:twist}) above
hold more generally
when $X$ is an arbitrary normal variety
of finite type over a finite field
of characteristic~$p$.
The proof is the same,
using \cite{Lafforgue-LanglandsCorresp}, Prop.~VII.7
(see~(\ref{rmk:generalization for 2abc}))
instead of~(\ref{thm:Lafforgue 12abc})
in the final step.
\end{remark}

\begin{proposition}
\label{prop:algsheaf}
Let $X$ be a smooth curve
over a finite field $\bF$ of characteristic~$p$.
Let $\ell \ne p$ be a prime number,
and let $\Lcal$ be a lisse $\bQellbar$-sheaf on $X$.
Suppose $\Lcal$ is algebraic.
Then $\Lcal$ is:
\begin{itemize}
\item[--]
$E$-rational, for some number field $E \subset\bQellbar$,
\item[--]
$C$-bounded in valuation, for some real number $C>0$, and
\item[--]
$D$-bounded in denominator, for some integer $D>0$.
\end{itemize}
\end{proposition}

\begin{proof}
It suffices to treat the case
when $\Lcal$ is irreducible,
and in that case,
the proposition
follows from the equivalence
$\text{(\ref{item:twist(a)})}\Leftrightarrow
 \text{(\ref{item:twist(f)})}$
in~(\ref{prop:twist}),
together with
Lafforgue's theorem~(\ref{thm:Lafforgue 12abc})
(assertions~(1) and~(2c))
and~(\ref{thm:2d})
(assertion~(2d)).
\end{proof}

\begin{theorem}
\label{thm:indeplalg}
Let $X$ be a smooth curve
over a finite field $\bF$ of characteristic~$p$.
Let $\ell \ne p$ be a prime number,
and let $\Lcal$ be
an irreducible (resp.~semisimple)
lisse $\bQellbar$-sheaf on $X$.
Suppose $\Lcal$ is
plain of characteristic~$p$.
Then
there exists a number field $E \subset \bQellbar$
and an absolutely irreducible (resp.~semisimple)
$E$-compatible system $\Lbold = \{\Lcal_\lambda\}$ on $X$
which is $E$-compatible
with the given lisse $\bQellbar$-sheaf $\Lcal$.
If $\Lcal$ is pure of weight~$w$ for some $w\in\bZ$,
then so is the system $\Lbold$.
\end{theorem}

\begin{proof}
It is clear that
we only have to treat the case
when $\Lcal$ is irreducible.
By the equivalence
$\text{(\ref{item:twist(a)})}\Leftrightarrow
 \text{(\ref{item:twist(f)})}$
in~(\ref{prop:twist}),
$\Lcal$ is isomorphic to a lisse $\bQellbar$-sheaf
of the form
$\Fcal\otimes(1/\alpha)^{\deg_{\bF}}$
where
$\Fcal$ is an irreducible lisse $\bQellbar$-sheaf
whose determinant is of finite order,
and $\alpha\in\bQellbar^\times$
is plain of characteristic~$p$.
Applying~(\ref{thm:Lafforgue 3}) to $\Fcal$,
we obtain a number field $E$ with $\alpha\in E$,
and an $E$-compatible system
$\Fbold = \{\Fcal_\lambda\}$
which is $E$-compatible with $\Fcal$.
The $E$-compatible system
$\Lbold = \{\Lcal_\lambda\}$,
with
$\Lcal_\lambda := \Fcal_\lambda\otimes(1/\alpha)^{\deg_{\bF}}$
for each $\lambda\in|E|_{\ne p}$,
satisfies the first assertion of the theorem.
The second assertion is clear.
\end{proof}

Note that
by the equivalence
$\text{(\ref{item:twist(a)})}\Leftrightarrow
 \text{(\ref{item:twist(b)})}$
in~(\ref{prop:twist}),
in order to verify the hypotheses that
$\Lcal$ is plain of characteristic~$p$
(resp.~is pure of weight $w$)
in~(\ref{thm:indeplalg}),
it suffices to
pick one closed point $x\in|X|$ of $X$
and check that
the $\bQellbar$-representation $[\Lcal(x)]$
of $\Gal(\kxbar/\kx)$
has the corresponding property.
Likewise,
in~(\ref{prop:algsheaf}),
it suffices to
check that
$[\Lcal(x)]$ is algebraic
for one closed point $x$
in order to conclude that
$\Lcal$ is algebraic.


\section{Frobenius tori}
\label{sect:Frobenius Tori}

In this section,
we review the theory of Frobenius tori
due to J.-P.~Serre,
with the aim of establishing
the two key results we need
---
namely,
(\ref{thm:Serre max torus})
and~(\ref{cor:max dim torus}).
The main ideas
are all in \cite{Serre-LettresRibet},~\S\S4--5;
we only have to
make them applicable to
general lisse $\ell$-adic sheaves
satisfying appropriate hypotheses.

\subsection{}
\label{sect:Frob tori}
Let $\Lambda$ be a field of characteristic~$0$,
and let $E \subset \Lambda$ be
a subfield of $\Lambda$.
Consider
an $r$-dimensional
$E$-rational $\Lambda$-representation $\sigma$
of $\Gal(\bFbar/\bF)$,
where $\bF$ denotes
an arbitrary finite field.
The image $\sigma(\Frob)\in\GL_r(\Lambda)$
of the geometric Frobenius
has a semisimple part
$$ \sigma(\Frob)^\semisimp \in \GL_r(\Lambda).
$$
Since $\sigma$ is $E$-rational,
the $E$-algebra generated by
the semisimple element $\sigma(\Frob)^\semisimp$,
$$ E[\sigma(\Frob)^\semisimp]  \ :=\ 
   E[T] \ / \  \left(\,
	       \text{minimal $E$-polynomial of $\sigma(\Frob)^\semisimp$
	      	     in the variable $T$} \,\right) ,
$$
is a finite etale algebra over $E$.
Therefore,
the multiplicative group of this $E$-algebra
defines an $E$-torus $\Mult(E[\sigma(\Frob)^\semisimp])$.
The group of $E$-points of $\Mult(E[\sigma(\Frob)^\semisimp])$
is isomorphic to $E[\sigma(\Frob)^\semisimp]^\times$,
and so it contains the element $\sigma(\Frob)^\semisimp$.
We define the {\sl Frobenius group $H(\sigma)$ of $\sigma$\/}
as
the Zariski closure
of the subgroup of $\Mult(E[\sigma(\Frob)^\semisimp])$
generated by 
$\sigma(\Frob)^\semisimp$.
It is a diagonalizable $E$-group,
whose group of characters $\tX(H(\sigma))$
as a $\Gal(\Ebar/E)$-module
is canonically isomorphic
---
via the ``evaluation at $\Frob$'' map
---
to the subgroup of $\Ebar^\times$
generated by the eigenvalues of $\sigma(\Frob)$.

One can also consider
the Zariski closure
of the subgroup of $\GL_r(\Lambda)$
generated by the element
$\sigma(\Frob)^\semisimp$:
$$ \overline{(\sigma(\Frob)^\semisimp)^\bZ}
   \quad\subseteq\GL_r(\Lambda);
$$
it is canonically isomorphic to
the $\Lambda$-scalar extension $H(\sigma)_{/ \Lambda}$
of the $E$-group $H(\sigma)$.

The {\sl Frobenius torus $A(\sigma)$ of $\sigma$\/}
is by definition
the identity component $H(\sigma)^0$ of $H(\sigma)$;
it is an $E$-torus,
whose group of characters $\tX(A(\sigma))$
as a $\Gal(\Ebar/E)$-module
is canonically isomorphic to
the torsion-free quotient of $\tX(H(\sigma))$.

Note that our definition
differs from that in \cite{Serre-LettresRibet}
in that
we first pass to the semi\-simplification
of the element $\sigma(\Frob)$;
this allows our discussion to proceed
without having to assume that
$\sigma$ is a semisimple representation.

\subsection{}
\label{sect:compat Frob tori}
The Frobenius group $H(\sigma)$
and the Frobenius torus $A(\sigma)$
of the $E$-rational representation $\sigma$
are determined
by the minimal $E$-polynomial
of the semisimple element $\sigma(\Frob)^\semisimp$,
and hence also
by the ``inverse'' characteristic polynomial
$$ \det(1-T\,\sigma(\Frob)) \in E[T]
$$
of $\sigma(\Frob)$.
Consequently,
if $\Lambda'$ is
another field containing $E$
and $\sigma'$ is
a $\Lambda'$-representation
of $\Gal(\bFbar/\bF)$
which is $E$-compatible with $\sigma$,
then the Frobenius groups $H(\sigma)$ and $H(\sigma')$
of $\sigma$ and $\sigma'$ respectively
are in fact the same $E$-group;
likewise
the Frobenius tori $A(\sigma)$ and $A(\sigma')$
are the same $E$-tori.

\subsection{}
\label{sect:order Frob finite quot}
The order of the finite group
$H(\sigma) / A(\sigma)$
is an integer which divides
the order of some root of unity
belonging to the field $E'$
obtained from $E$ by
adjoining all eigenvalues of $\sigma(\Frob)$,
and $E'$ is of degree~$\le r!$ over $E$.
Suppose $E$ is a number field;
then there are only finitely many
roots of unity of degree~$\le r!$ over $E$.
Letting $N$ be the least common multiple
of the orders of these roots of unity,
we see that
the order of $H(\sigma) / A(\sigma)$
divides $N$.
Note that
$N$ depends only on
the number field $E$ and the integer~$r>0$,
but is independent of $\sigma$.

\subsection{}
In order to state
the finiteness result of Serre,
let us introduce the following notation.
Let $r\ge 1$ be an integer,
and let $\Lambda$ be a field of characteristic~$0$.
Let $E \subset \Lambda$ be
a number field contained in $\Lambda$.
Let $C \ge 0$ be a real number,
and let $D > 0$ be an integer.
Define the set
$$ {\Sigma}^r_{E\subset\Lambda}(C,D)_p :=
 \left\{\quad \text{\parbox{26em}{
$r$-dimensional $\Lambda$-representations
$\sigma$ of $\Gal(\bFbar/\bF)$
such that
for some positive power $q$ of $p$,
and some integer $w$,
$\sigma$ is:
\begin{enumerate}
\item[(1)]
$E$-rational,
\item[(2a)]
pure of weight~$w$ with respect to~$q$,
\item[(2b)]
plain of characteristic~$p$,
\item[(2c)]
$C$-bounded in valuation with respect to~$q$, and
\item[(2d)]
$D$-bounded in denominator with respect to~$q$
\end{enumerate}
   }} \quad\right\}
$$
and let
$$ {\Sigma}^r_{E\subset\Lambda}(C,D) :=
   \bigcup_{\text{$p$ prime}}  {\Sigma}^r_{E\subset\Lambda}(C,D)_p.
$$

\begin{theorem}[J.-P.~Serre]
\label{thm:Serre finiteness}
With the above notation,
as $\sigma$ runs over
the set $\Sigma^r_{E\subset\Lambda}(C,D)$,
there are only finitely many possibilities for
the ${\GL_r}_{/ \Ebar}$-conjugacy class
of the Frobenius torus $A(\sigma)$ of $\sigma$,
where $\Ebar$ denotes
an algebraic closure of $E$.
\end{theorem}

We leave to the reader
the pleasant exercise of generalizing
the argument found in~\cite{Serre-LettresRibet}, Th.~on p.~8
to the situation of~(\ref{thm:Serre finiteness});
the case of weight $w=0$
requires a little more argument.
We note in passing that
the assumption that
$C$ and $D$ are {\it uniform bounds\/}
for the valuation and denominator
of the $\sigma$ in question
plays a key role
in showing the finiteness assertion here.

\subsection{}
\label{sect:Frob torus notation}
Now let $X$ be an irreducible normal scheme
of finite type over
a finite field of characteristic~$p$,
and let $\etabar \rightarrow X$ be
a geometric point of $X$.
Assume that
$\Lambda$ is
an $\ell$-adic field
for some prime number~$\ell\ne p$.
Let $\Lcal$ be
a {\it semisimple\/} lisse $\Lambda$-sheaf on $X$;
its arithmetic monodromy group
$\Garith(\Lcal,\etabar)$ of $\Lcal$
is then
a (possibly non-connected) reductive group over $\Lambda$.
For every closed point $x\in|X|$,
the $\Lambda$-representation $[\Lcal(x)]$
of $\Gal(\kxbar/\kx)$
factors through $\Garith(\Lcal,\etabar)$;
the semisimple element
$[\Lcal(x)](\Frob_x)^\semisimp$
therefore lies in $\Garith(\Lcal,\etabar)$.
Hence, for each $x\in|X|$,
we have the diagonalizable $\Lambda$-group
\begin{equation*}
 H(x) := \overline{([\Lcal(x)](\Frob_x)^\semisimp)^\bZ}
 \quad\subseteq\Garith(\Lcal,\etabar)
\end{equation*}
and the $\Lambda$-torus
\begin{equation*}
 A(x) := H(x)^0
 \quad\subseteq\Garith(\Lcal,\etabar)^0.
\end{equation*}
These are canonically isomorphic, over $\Lambda$,
to the Frobenius group $H([\Lcal(x)])_{/ \Lambda}$
and the Frobenius torus $A([\Lcal(x)])_{/ \Lambda}$
associated to the representation $[\Lcal(x)]$
of $\Gal(\kxbar/\kx)$.

\begin{theorem}
\label{thm:Serre max torus}
{\rm (Cf.~\cite{Serre-LettresRibet}, Th.~on p.~12.)}
Assume the notation and hypotheses
of~(\ref{sect:Frob torus notation}),
and let $r \ge 1$ be
the rank of
the semisimple lisse $\Lambda$-sheaf $\Lcal$ 
on $X$.
Suppose there exist
a number field $E \subset \Lambda$
contained in the $\ell$-adic field $\Lambda$,
an integer~$w$,
a real number~$C \ge 0$,
and an integer~$D > 0$,
such that
the lisse $\Lambda$-sheaf $\Lcal$ on $X$ is
\begin{enumerate}
\item[(1)]
$E$-rational,
\item[(2a)]
pure of weight~$w$,
\item[(2b)]
plain of characteristic~$p$,
\item[(2c)]
$C$-bounded in valuation, and
\item[(2d)]
$D$-bounded in denominator.
\end{enumerate}
Then
there exists
an open dense subset
$U \subset \Garith(\Lcal,\etabar)^0$
of the identity component of
$\Garith(\Lcal,\etabar)$,
such that
$U$ is stable under
conjugation by $\Garith(\Lcal,\etabar)$,
and such that
for any closed point $x \in |X|$ of $X$
with $[\Lcal(x)](\Frob_x) \in U$,
the Zariski closure 
$$ \overline{([\Lcal(x)](\Frob_x))^\bZ}
 \subseteq \Garith(\Lcal,\etabar)^0
$$
of the subgroup it generates
is a maximal $\Lambda$-torus
of $\Garith(\Lcal,\etabar)^0$.
\end{theorem}

\begin{proof}
We may assume that
$\Lambda$ is algebraically closed.
By the list of hypotheses on
the lisse $\Lambda$-sheaf $\Lcal$,
we may apply~(\ref{thm:Serre finiteness}) to
the collection of repre\-sentations $[\Lcal(x)]$
($x\in|X|$),
and it follows that
as $x$ runs over $|X|$,
there are only finitely many possibilities for
the ${\GL_r}$-conjugacy class
of the Frobenius torus $A([\Lcal(x)])$
associated to $[\Lcal(x)]$.

Let $T_0 \subseteq \Garith(\Lcal,\etabar)^0$
be a maximal torus of $\Garith(\Lcal,\etabar)^0$.
We claim that
each ${\GL_r}$-conjugacy class
of Frobenius tori
contains only finitely many
subtori $A \subseteq T_0$ of $T_0$.
To see this,
we just have to show that
if two subtori $A_1, A_2$ of $T_0$
are conjugate to each other in $\GL_r$,
then they are already conjugate to each other
under the action of
the (finite) Weyl group of $T_0$ in ${\GL_r}$.
But if $g\in\GL_r$ conjugates
$A_1\subseteq T_0$ to $A_2\subseteq T_0$,
then both $T_0$ and the $g$-conjugate of $T_0$
are maximal tori contained in
the centralizer of $A_2$ in $\GL_r$;
this centralizer
is a connected reductive subgroup of $\GL_r$,
so one can adjust $g$
by an element of this subgroup
to assume that
$g$ normalizes $T_0$,
and this shows what we want.
As a consequence of our claim,
we see that
if we define the set
$$ \Phi' :=
   \left\{\quad \text{\parbox{25em}{
proper subtori $A \subsetneq T_0$ of $T_0$
such that, for some $x\in|X|$,
$A$ belongs to the ${\GL_r}$-conjugacy class
of the Frobenius torus $A([\Lcal(x)])$
associated to $[\Lcal(x)]$
   }} \quad\right\},
$$
then $\Phi'$ is a finite set.

We are given the number field $E$
and the integer~$r>0$.
Let $N$ be the least common multiple of
the orders of
the finitely many roots of unity
of degree~$\le r!$ over $E$.
Define the set
$$ \Phi :=
   \left\{\quad \text{\parbox{19em}{
subgroups $H \subseteq T_0$ of $T_0$
such that\\
the identity component $H^0$ of $H$
lies in $\Phi'$\\
and the order of $H/H^0$ divides $N$
   }}  \quad\right\}.
$$
From the finiteness of the set $\Phi'$,
we see that 
$\Phi$ is also a finite set.

If $H \subseteq T_0$
belongs to the finite set $\Phi$,
then its identity component $H^0$
is a proper subgroup of $T_0$,
and so the same is true for $H$.
Therefore, if $F_H$ denotes
the Zariski closure of the union of
all $\Garith(\Lcal,\etabar)$-conjugates of $H$,
then $F_H$ is of dimension
strictly less than
that of $\Garith(\Lcal,\etabar)^0$.
It follows that
$$ U' := \Garith(\Lcal,\etabar)^0
	- \bigcup_{H \in \Phi} F_H
$$
is an open dense subset
of $\Garith(\Lcal,\etabar)^0$.
We define $U$ to be
the intersection of $U'$
with the open dense (``characteristic'') subset of
regular semisimple elements in
the connected reductive group $\Garith(\Lcal,\etabar)^0$.
From its construction
we see that
$U$ is stable under
conjugation by $\Garith(\Lcal,\etabar)$.

Let $x\in|X|$ be a closed point
with $[\Lcal(x)](\Frob_x) \in U$.
Since $[\Lcal(x)](\Frob_x)$ is semisimple,
the Zariski closure 
$$ \overline{([\Lcal(x)](\Frob_x))^\bZ}
	 \quad\subseteq \Garith(\Lcal,\etabar)^0
$$
of the subgroup generated by $[\Lcal(x)](\Frob_x)$
is just
the diagonalizable $\Lambda$-group $H(x)$
introduced in~(\ref{sect:Frob torus notation}),
and it is contained in some maximal torus
of $\Garith(\Lcal,\etabar)^0$.
The identity component
of $H(x)$
is the torus $A(x) = H(x)^0$,
and by construction,
it belongs to
the ${\GL_r}$-conjugacy class
of the Frobenius torus $A([\Lcal(x)])$
associated to $[\Lcal(x)]$.
From the assumption that
the polynomial
$$ \det(1-T\,[\Lcal(x)](\Frob_x))
$$
has coefficients in the number field $E$
and has degree $r$,
we see (cf.~(\ref{sect:order Frob finite quot})) that
the order of $H(x) / A(x)$
divides the integer~$N$ above.
If $A(x)$ is 
not a maximal torus of $\Garith(\Lcal,\etabar)^0$,
then
$H(x)$ is $\Garith(\Lcal,\etabar)^0$-conjugate
to some subgroup $H$
belonging to the set $\Phi$,
but this contradicts the fact that
$[\Lcal(x)](\Frob_x)$
does not lie in $F_H$.
Thus $A(x)$ must be a maximal torus
of $\Garith(\Lcal,\etabar)^0$,
and hence
the same is true for $H(x)$.
\end{proof}

\begin{corollary}
\label{cor:max dim torus}
Assume the notation and hypotheses
of~(\ref{sect:Frob torus notation})
and~(\ref{thm:Serre max torus}).
Let $x_0\in|X|$ be a closed point of $X$.
The subgroup $A(x_0) \subseteq\Garith(\Lcal,\etabar)$
is a maximal torus of $\Garith(\Lcal,\etabar)^0$
if and only if
the Frobenius torus $A([\Lcal(x_0)])$
associated to $[\Lcal(x_0)]$
has the maximum dimension (over $E$)
among all Frobenius tori $A([\Lcal(x)])$
as $x$ runs over the set of closed points $|X|$ of $X$.
\end{corollary}

\begin{proof}
By~(\ref{thm:Serre max torus})
and \u{C}ebotarev's density theorem
(cf.~\cite{Serre-ZetaLFunctions}, Th.~7),
there exists some closed point $y\in|X|$
such that
$A(y)$ is a maximal torus
of $\Garith(\Lcal,\etabar)^0$;
so the dimension of each $A(x)$ (over $\Lambda$)
is at most that of $A(y)$.
Therefore,
$A({x_0})$ is a maximal torus
of $\Garith(\Lcal,\etabar)^0$
if and only if
it has the maximum dimension (over $\Lambda$)
among the tori $A(x)$
as $x$ runs over $|X|$.
Since one has the isomorphism
$A(x) \cong A([\Lcal(x)])_{/ \Lambda}$
for each $x\in|X|$,
the corollary follows.
\end{proof}

\begin{remark}
The results of~(\ref{thm:Serre max torus})
and~(\ref{cor:max dim torus})
also hold
when $X$ is any
irreducible normal scheme of finite type
over $\bZ[1/\ell]$,
except that one would have to replace
the assumption~(2b)
on the lisse sheaf $\Lcal$
by an appropriately generalized
notion of plainness
(one in which the characteristic $p$
 is allowed to vary).
The proof is the same,
thanks to the very general nature
of~(\ref{thm:Serre finiteness}).
\end{remark}


\section{Independence of~$\ell$}
\label{sect:Indep of ell}

In this section,
we prove
our main ``independence of~$\ell$'' theorem~(\ref{thm:main arith}).

\subsection{Hypotheses}
\label{sect:hypotheses}
Let $X$ be a smooth curve
over a finite field of characteristic~$p$,
let $E$ be a number field,
and let $\Lbold = \{\Lcal_\lambda\}$ be
an $E$-compatible system of lisse sheaves
on the curve $X$.
We assume that
the $E$-compatible system $\Lbold$ is
semisimple and
pure of weight~$w$ for some integer~$w$.
For each $\lambda\in|E|_{\ne p}$,
write
$$ G_{E_\lambda} := \Garith(\Lcal_\lambda,\etabar)
$$
for the arithmetic monodromy group of $\Lcal_\lambda$,
and write $\sigma_\lambda$
for its tautological faithful representation.
The group $G_{E_\lambda}$
is a (possibly non-connected)
reductive group over $E_\lambda$,
and the representation $\sigma_\lambda$
is an $E_\lambda$-rational representation of $G_{E_\lambda}$.

Our goal is to show that
the identity component $G_{E_\lambda}^0$ of $G_{E_\lambda}$
and its tautological representation $\sigma_\lambda$
are ``independent of $\lambda$''
in the sense of~(\ref{thm:main arith}).
In~(\ref{sect:constructions}),
we shall construct a finite extension $F$ of $E$,
a connected split reductive group $G_0$ over $F$,
and an $F$-rational representation $\sigma_0$ of $G_0$.
We will show in~(\ref{cor:stronger main thm arith}) that
for every place $\lambda\in|F|_{\ne p}$ of $F$
not lying over~$p$,
writing $\lambda$ also for its restriction to $E$,
one has an isomorphism
of connected $F_\lambda$-groups:
$$ f_\lambda:
   G_{E_\lambda}^0 \otimes_{E_\lambda} F_\lambda
   \xrightarrow{\cong}
   G_0 \otimes_F F_\lambda,
$$
and we will show in~(\ref{proof:main(ii)}) that
with this isomorphism of groups,
one also has an isomorphism
of representations:
$$ \sigma_\lambda \otimes_{E_\lambda} F_\lambda
   \cong
   \sigma_0 \otimes_F F_\lambda.
$$

First, some preliminary reductions.

\begin{lemma}
\label{lemma:reduction to ctd}
To prove~(\ref{thm:main arith}),
we may assume that
for each $\lambda\in|E|_{\ne p}$,
the $E_\lambda$-monodromy group $G_{E_\lambda}$
is {\rm connected}.
\end{lemma}

\begin{proof}
Indeed, under the hypotheses of~(\ref{sect:hypotheses}),
we may apply Serre's theorem~(\ref{thm:Serre indep kernel arith})
to see that
as $\lambda$ runs over $|E|_{\ne p}$,
the kernel of the surjective homomorphism
$$ \pi_1(X,\etabar)
   \xrightarrow{[\Lcal_\lambda]}
   G_{E_\lambda}
   \rightarrow\joinrel\rarrow
   G_{E_\lambda} / G_{E_\lambda}^0
$$
is the same open subgroup of $\pi_1(X,\etabar)$,
independent of $\lambda\in|E|_{\ne p}$.
This open subgroup of $\pi_1(X,\etabar)$
corresponds to
a finite etale cover $a: X' \rightarrow X$ of $X$
by a connected smooth curve $X'$
pointed by the same geometric point $\etabar$.
The inverse image $a^*\Lbold := \{a^*\Lcal_\lambda\}$ on $X'$
of the $E$-compatible system $\Lbold = \{\Lcal_\lambda\}$
is still a semisimple $E$-compatible system
which is pure of integer weight.
Moreover, for each $\lambda\in|E|_{\ne p}$,
the arithmetic monodromy group of $a^*\Lcal_\lambda$
is the connected reductive group $G_{E_\lambda}^0$,
given with the same
tautological representation $\sigma_\lambda$.
Hence, to prove~(\ref{thm:main arith}),
we may replace $X$ by $X'$
and $\Lbold$ by its inverse image.
This proves the lemma.
\end{proof}

\subsection{Notation}
\label{sect:notation}
Henceforth,
we adopt
both the hypotheses~(\ref{sect:hypotheses})
and the assumption that
for each $\lambda\in|E|_{\ne p}$,
we have $G_{E_\lambda} = G_{E_\lambda}^0$.

Let $x\in|X|$ be a closed point of $X$.
By the $E$-compatibility hypothesis on $\Lbold$,
as $\lambda$ runs over $|E|_{\ne p}$,
the Frobenius group $H([\Lcal_\lambda(x)])$
(cf.~(\ref{sect:Frob tori}))
associated to
the $E_\lambda$-representation $[\Lcal_\lambda(x)]$
of $\Gal(\kxbar/\kx)$
is a diagonalizable group over $E$
which is
independent of the place $\lambda\in|E|_{\ne p}$
(cf.~(\ref{sect:compat Frob tori}));
we denote it by $H(\Lbold(x))$.
Similarly,
the Frobenius torus $A([\Lcal_\lambda(x)])$
is an $E$-torus
which is
independent of the place $\lambda\in|E|_{\ne p}$
and we denote it by $A(\Lbold(x))$;
it is also the identity component of $H(\Lbold(x))$.

For each $\lambda\in|E|_{\ne p}$,
we write
$$ H(x)_{E_\lambda} :=
   \overline{([\Lcal_\lambda(x)](\Frob_x)^\semisimp)^\bZ}
   \quad\subseteq G_{E_\lambda}
$$
for the Zariski closure
of the subgroup of $G_{E_\lambda}$
generated by $[\Lcal(x)](\Frob_x)^\semisimp$.
It is
a diagonalizable subgroup of $G_{E_\lambda}$
over $E_\lambda$,
and we write
$$ A(x)_{E_\lambda} := H(x)_{E_\lambda}^0
   \quad\subseteq G_{E_\lambda}
$$
for its identity component,
which is an $E_\lambda$-torus.

\begin{lemma}
\label{lemma:Frob max torus}
Assume~(\ref{sect:hypotheses}) and~(\ref{sect:notation}) above.
There exist
infinitely many closed points \linebreak
$x\in|X|$ of $X$
with the following property:
for every place $\lambda\in|E|_{\ne p}$ of $E$
not lying over~$p$,
one has
$ A(x)_{E_\lambda} = H(x)_{E_\lambda}$,
and it is
a maximal $E_\lambda$-torus of $G_{E_\lambda}$.
\end{lemma}

\begin{proof}
By our assumptions
on the $E$-compatible system $\Lbold = \{\Lcal_\lambda\}$,
for each 
$\lambda\in|E|_{\ne p}$,
the lisse $E_\lambda$-sheaf $\Lcal_\lambda$
is $E$-rational (hence algebraic),
pure of some integer weight,
and plain of characteristic~$p$;
by~(\ref{prop:algsheaf}), it is also
$C$-bounded in valuation
and
$D$-bounded in denominator.
Hence
the hypotheses of~(\ref{thm:Serre max torus})
and~(\ref{cor:max dim torus})
are all verified for these lisse sheaves.

Pick any place $\mu\in|E|_{\ne p}$
of $E$ not lying over~$p$.
From~(\ref{thm:Serre max torus}),
we obtain 
an open dense subset
$U \subset G_{E_\mu}$,
stable under conjugation by $G_{E_\mu}$,
such that
for any $x\in|X|$
with $[\Lcal_\mu(x)](\Frob_x)\in U$,
the subgroup $H(x)_{E_\mu}\subseteq G_{E_\mu}$
is a maximal $E_\mu$-torus of $G_{E_\mu}$
(recall from~(\ref{lemma:reduction to ctd}) that
we have put ourselves in the situation where
the monodromy group $G_{E_\lambda}$ is connected
for every $\lambda\in|E|_{\ne p}$).
By \u{C}ebotarev's density theorem,
there exist infinitely many
closed points $x\in|X|$ of $X$
with $[\Lcal_\mu(x)](\Frob_x)\in U$.
We claim that
these closed points
satisfy the conclusion of the lemma.

To see that,
let $x_0$ be such a closed point.
By~(\ref{cor:max dim torus}) and the assumption on $x_0$,
we infer that
$A([\Lbold(x_0)]) = H([\Lbold(x_0)])$
and that
the Frobenius torus $A([\Lbold(x_0)])$
has the maximum dimension (over $E$)
among all Frobenius tori $A([\Lbold(x)])$
as $x$ runs over $|X|$.
Thus, for any place $\lambda\in|E|_{\ne p}$ of $E$
not lying over~$p$,
it follows from~(\ref{cor:max dim torus}) again that
the $E_\lambda$-subgroup
$A(x_0)_{E_\lambda} = H(x_0)_{E_\lambda} \subseteq G_{E_\lambda}$
is a maximal torus of $G_{E_\lambda}$.
\end{proof}

\subsection{Constructions}
\label{sect:constructions}
Assume the hypotheses in~(\ref{sect:hypotheses})
and the notation in~(\ref{sect:notation}).
We shall now make a series of constructions
to be used for the rest of
this section.

\subsubsection{}
\label{subsect:closed point}
Choose once and for all
a closed point $x\in|X|$ of $X$
satisfying the conclusion of~(\ref{lemma:Frob max torus});
thus for every place $\lambda\in|E|_{\ne p}$
of $E$ not lying over~$p$,
$A(x)_{E_\lambda} = H(x)_{E_\lambda}$
is a maximal torus of $G_{E_\lambda}$.

\subsubsection{}
\label{subsect:splitting field}
Let $F$ be the splitting field of
the $E$-torus $A([\Lbold(x)]) = H([\Lbold(x)])$;
it is the same as the splitting field of the polynomial
$\det(1-T\,\Frob_x,\Lcal_\lambda) \in E[T]$, 
for any/every $\lambda\in|E|_{\ne p}$.

\subsubsection{}
\label{subsect:split torus}
Let $T_0$ be the split $F$-torus
defined by the Frobenius $E$-torus $A([\Lbold(x)])$
via scalar extension:
$$ T_0 := A([\Lbold(x)]) \otimes_E F
        = H([\Lbold(x)]) \otimes_E F.
$$

\subsubsection{}
\label{subsect:scalar extn}
For each $\lambda\in|F|_{\ne p}$,
writing $\lambda$ also for its restriction to $E$,
we obtain the lisse $F_\lambda$-sheaf
$\Fcal_\lambda := \Lcal_\lambda \otimes_{E_\lambda} F_\lambda$
by scalar extension.
This gives us
the collection of lisse sheaves
$\Fbold := \{\Fcal_\lambda\}$
indexed by the places $|F|_{\ne p}$ of $F$
not lying over~$p$,
and it is clear that
$\Fbold$ is an $F$-compatible system.

\subsubsection{}
\label{subsect:monodromy group}
For each $\lambda\in|F|_{\ne p}$,
the arithmetic monodromy group
$\Garith(\Fcal_\lambda,\etabar)$
of the lisse $F_\lambda$-sheaf $\Fcal_\lambda$
is identified with
$G_{E_\lambda} \otimes_{E_\lambda} F_\lambda$;
we denote it by $G_{F_\lambda}$.
Thus:
$$ G_{F_\lambda}
 := G_{E_\lambda} \otimes_{E_\lambda} F_\lambda
 =  \Garith(\Fcal_\lambda,\etabar).
$$
The tautological representation
of $G_{F_\lambda}$
is identified with
$\sigma_\lambda\otimes_{E_\lambda} F_\lambda$;
we denote this by $\rho_\lambda$.
Thus:
$$ \rho_\lambda :=
   \sigma_\lambda \otimes_{E_\lambda} F_\lambda.
$$

\subsubsection{}
\label{subsect:max torus}
For each $\lambda\in|F|_{\ne p}$,
writing $\lambda$ also for its restriction to $E$,
we have the $F_\lambda$-torus
$$ T_{F_\lambda}
 := A(x)_{E_\lambda} \otimes_{E_\lambda} F_\lambda
  = H(x)_{E_\lambda} \otimes_{E_\lambda} F_\lambda
 \quad\subseteq G_{F_\lambda}.
$$
By~(\ref{subsect:closed point}) and~(\ref{subsect:splitting field}),
$T_{F_\lambda}$
is in fact a {\it split} maximal torus
of $G_{F_\lambda}$,
whence $G_{F_\lambda}$
is a connected {\it split} reductive group
over $F_\lambda$.

\subsubsection{}
\label{subsect:torus isom}
For each $\lambda\in|F|_{\ne p}$,
one has the canonical isomorphism (cf.~(\ref{sect:Frob tori}))
$$ f_\lambda:
   T_{F_\lambda}
   \xrightarrow{\cong}
   {T_0}_{/ F_\lambda}
 \qquad\text{of split $F_\lambda$-tori},
$$
where the right-hand side 
${T_0}_{/ F_\lambda}$
is the $F_\lambda$-scalar extension
of the $F$-torus $T_0$.

\subsubsection{}
\label{subsect:mu}
Choose once and for all
a place $\mu\in|F|_{\ne p}$ of $F$
not lying over~$p$.

\subsubsection{}
\label{subsect:G0}
Define
the $F$-group $G_0$
as follows.
First consider
the $F_\mu$-group
$G_{F_\mu}$;
it is a connected split reductive group
over $F_\mu$,
and so it is defined over
any subfield of $F_\mu$.
We let $G_0$ be
``the'' connected split reductive group over $F$,
whose isomorphism type over $F_\mu$
is that of $G_{F_\mu}$,
and which contains
the split $F$-torus $T_0$
as its maximal torus.
Thus, $G_0$ and $T_0$
fit into the commutative diagram
$$\begin{matrix}
  T_{F_\mu} &
  \rightarrow &
  T_0 \\
 \hookdownarrow & &
 \hookdownarrow \\
  G_{F_\mu} &
  \rightarrow &
  G_0 \\
  \downarrow & &
  \downarrow \\
  \Spec(F_\mu) &
  \rightarrow &
  \Spec(F) 
\end{matrix}
$$
in which the squares are cartesian.

\subsubsection{}
\label{subsect:sigma0}
Likewise,
define the $F$-rational representation $\sigma_0$
of the $F$-group $G_0$ as follows.
The $F_\mu$-group $G_{F_\mu}$
is given with
a tautological representation $\rho_\mu$,
and this is defined over
any subfield of $F_\mu$
because $G_{F_\mu}$ is split.
We let $\sigma_0$ be
``the'' $F$-rational representation
of the $F$-group $G_0$
whose scalar-extension to $F_\mu$
is the representation $\rho_\mu$.
Thus:
$$ \rho_\mu :=
   \sigma_\mu \otimes_{E_\mu} F_\mu
   \cong
   \sigma_0 \otimes_F F_\mu.
$$

Our goal now is
to show that
for each $\lambda\in|F|_{\ne p}$,
the canonical isomorphism $f_\lambda$
in~(\ref{subsect:torus isom})
of split $F_\lambda$-tori
{\sl extends\/} to
an isomorphism
$$ f_\lambda:
   G_{F_\lambda}
   \xrightarrow{\cong}
   {G_0}_{/ F_\lambda}
\qquad\text{of connected split reductive groups over $F_\lambda$}
$$
(where the right-hand side 
${G_0}_{/ F_\lambda}$
is the $F_\lambda$-scalar extension
of the $F$-group $G_0$),
and  that
if we identify the two groups
using this isomorphism,
then we have an isomorphism of representations
$$ \rho_\lambda :=
   \sigma_\lambda \otimes_{E_\lambda} F_\lambda
   \cong
   \sigma_0 \otimes_F F_\lambda.
$$
These will be achieved
in~(\ref{cor:stronger main thm arith})
and in~(\ref{proof:main(ii)})
respectively.

\subsection{Notation}
\label{sect:more notation}
If $G$ is
a split reductive algebraic group
over a field $k$,
and $T \subseteq G$ is
a maximal torus of $G$,
we denote
the group of characters of $T$ over $k$
by $\tX_k(T)$,
and we set
$$ \Irr_k(G) :=
 \left\{\quad
 \text{isomorphism classes of
       irreducible $k$-rational representations of $G$}
 \quad\right\}.
$$
Note that
since $G$ is split,
every object of $\Irr_k(G)$
is in fact absolutely irreducible.
We write $\tK(T)$ and $\tK(G)$ for
the Grothendieck ring of
the category of $k$-rational representations of
$T$ and $G$ respectively;
these rings are
given with the canonical basis of
$\tX_k(T)$ and $\Irr_k(G)$ respectively.
Restriction of representations
from $G$ to $T$
gives rise to
the ``character homomorphism''
$\ch_G : \tK(G) \rightarrow \tK(T)$
of Grothendieck rings;
this homomorphism is injective
if $G$ is connected.

The main representation-theoretic input
we need is the fact that
over a field of characteristic~$0$,
the isomorphism type of
a connected split reductive group
can be determined
from the character data of its
irreducible representations;
more precisely:

\begin{theorem}
\label{thm:determine gp}
Let $k$ be a field of characteristic~$0$.
Let $G$ and $G'$ be
connected split
reductive algebraic groups over $k$,
and let $T \subseteq G$ and $T' \subseteq G'$ be
maximal tori,
satisfying the hypotheses~(i)--(iii) below.
\begin{itemize}
\item[(i)]
Suppose
$\phi: \Irr_k(G') \xrightarrow{\cong} \Irr_k(G)$
is a bijection of sets;
it induces
an isomorphism of abelian groups
$ \phi: \tK(G') \xrightarrow{\cong} \tK(G) $
which makes the diagram
$$
\begin{matrix}
   \tK(G')    &
   \xrightarrow{\overset{\ \ \phi\ \ }{\cong}} &
   \tK(G) \\
   \hookuparrow   & &
   \hookuparrow \\[-2ex]
   \Irr_k(G')  &
   \xrightarrow{\overset{\ \ \phi\ \ }{\cong}} &
   \Irr_k(G)
\end{matrix}
\qquad\text{commute}.
$$
\item[(ii)]
Suppose
$f_T : T \xrightarrow{\cong} T'$
is an isomorphism of tori;
it induces
an isomorphism of the group of characters
$ \tX(f_T) : \tX_k(T') \xrightarrow{\cong} \tX_k(T) $
and an isomorphism of rings
$ \tK(f_T) : \tK(T') \xrightarrow{\cong} \tK(T) $
which make the diagram
$$\begin{matrix}
   \tK(T')    &
   \xrightarrow{\overset{\tK(f_T)}{\cong}} &
   \tK(T) \\
   \hookuparrow & &
   \hookuparrow \\[-2ex]
   \tX_k(T')  &
   \xrightarrow{\overset{\tX(f_T)}{\cong}} &
   \tX_k(T)
\end{matrix}
\qquad\text{commute}.
$$
\item[(iii)]
Finally, suppose that
the diagram
$$\begin{matrix}
   \tK(T')    &
   \xrightarrow{\overset{\tK(f_T)}{\cong}} &
   \tK(T) \\
   \scriptstyle{\ch_{G'}} \ \hookuparrow & &
   \hookuparrow \ \scriptstyle{\ch_{G}}  \\[-2ex]
   \tK(G')    &
   \xrightarrow{\overset{\ \ \phi\ \ }{\cong}} &
   \tK(G)
\end{matrix}
\qquad\text{commutes}.
$$
\end{itemize}
Then the following conclusions hold:
\begin{enumerate}
\item
There exists
an isomorphism of algebraic groups
$ f: G \xrightarrow{\cong} G' $
such that:
\begin{enumerate}
\item
\label{output:isom}
$f$ extends the isomorphism $f_T$
of maximal tori
given in~(ii),
i.e.~the diagram
$$\begin{matrix}
  T  &
  \xrightarrow{\overset{f_T}{\cong}} &
  T' \\
 \hookdownarrow & &
 \hookdownarrow \\[-2ex]
  G  &
  \xrightarrow{\overset{f}{\cong}} &
  G'
\end{matrix}
\qquad\text{commutes}; \quad\text{and}
$$
\item
\label{output:Irr}
the bijection of sets
$$ \Irr_k(f) : \Irr_k(G') \xrightarrow{\cong} \Irr_k(G)
\qquad\text{given by $\rho' \mapsto \rho'\circ f$}
$$
is equal to
the bijection of sets
$$ \phi : \Irr_k(G') \xrightarrow{\cong} \Irr_k(G)
\qquad\text{given in~(i)}.
$$
\end{enumerate}
\item
\label{output:unique}
If $f' : G \xrightarrow{\cong} G'$
is another isomorphism of algebraic groups
having the same properties as $f$
in~(\ref{output:isom}) and~(\ref{output:Irr}) above,
then there exists a $k$-rational point $t \in T(k)$ of $T$
such that
$$ f' = f \circ \text{(conjugation by $t$)}.
$$
\end{enumerate}
\end{theorem}

This result is intuitively clear,
despite its lengthy statement.
It seems to be also well known to
the experts in algebraic groups,
although I have not been able to locate
a satisfactory reference in the published literature.
The prudent reader
can refer to \cite{Chin-DetermineRedGp},~Th.~1.4
for a detailed proof.

To proceed, 
let us assume for the moment
the validity of the following:

\begin{theorem}
\label{thm:bijection Irr}
Given the constructions~(\ref{sect:constructions}) above,
for each $\lambda\in|F|_{\ne p}$,
there exists a unique bijection
\begin{equation*}
\tag{\ref{thm:bijection Irr}.1}
\label{eqn:bijection}
 \phi:
   \Irr_{F_\lambda}({G_0}_{/ F_\lambda})
   \xrightarrow{\cong}
   \Irr_{F_\lambda}(G_{F_\lambda})
\qquad\text{of sets},
\end{equation*}
such that
the diagram
\begin{equation*}
\tag{\ref{thm:bijection Irr}.2}
\label{eqn:Irr}
\begin{matrix}
  \tK({T_0}_{/ F_\lambda}) &
  \xrightarrow{\overset{\tK(f_\lambda)}{\cong}} &
  \tK(T_{F_\lambda}) \\
  \scriptstyle{\ch_{{G_0}_{/ F_\lambda}}} \ \hookuparrow & &
  \hookuparrow \ \scriptstyle{\ch_{G_{F_\lambda}}} \\[-2ex]
  \tK({G_0}_{/ F_\lambda}) &
  \xrightarrow{\overset{\ \ \phi\ \ }{\cong}} &
  \tK(G_{F_\lambda})
\end{matrix}
\qquad\text{commutes}.
\end{equation*}
\end{theorem}

Now if we combine~(\ref{thm:determine gp})
with~(\ref{thm:bijection Irr}),
we obtain the following
strengthened version
of our main result~(\ref{thm:main arith})(i):

\begin{corollary}
\label{cor:stronger main thm arith}
Given the constructions~(\ref{sect:constructions}) above,
for each $\lambda\in|F|_{\ne p}$,
there exists
an isomorphism
$$ f_\lambda:
   G_{F_\lambda}
   \xrightarrow{\cong}
   {G_0}_{/ F_\lambda}
\qquad\text{of connected split reductive groups over $F_\lambda$},
$$
extending the canonical isomorphism
(cf.~(\ref{subsect:torus isom}))
$$  f_\lambda:
    T_{F_\lambda}
    \xrightarrow{\cong}
    {T_0}_{/ F_\lambda}
  \qquad\text{of split $F_\lambda$-tori},
$$
and inducing the bijection~(\ref{eqn:bijection})
given by~(\ref{thm:bijection Irr}).
Moreover, if
$$ f'_\lambda:
   G_{F_\lambda}
   \xrightarrow{\cong}
   {G_0}_{/ F_\lambda}
\qquad\text{is another isomorphism}
$$
with the same properties as $f_\lambda$ above,
then there exists $t\in T_{F_\lambda}(F_\lambda)$
such that
$$ f'_\lambda  =  f_\lambda \circ \text{(conjugation by $t$)}.
$$
\end{corollary}

We postpone the discussion of~(\ref{thm:main arith})(ii)
to~(\ref{proof:main(ii)});  
for now,
let us focus on
the proof of~(\ref{thm:bijection Irr}).

\subsection{}
\label{sect:reformulation}
Let $\lambda\in|F|_{\ne p}$ be
a given place of $F$
not lying over~$p$.
A bijection~(\ref{eqn:bijection})
which makes~(\ref{eqn:Irr}) commute
is necessarily unique,
as one sees from the fact that
the two sides of~(\ref{eqn:bijection})
are the canonical basis sets
for the bottom two terms of~(\ref{eqn:Irr}).
Hence it suffices for us to show
the existence of~(\ref{eqn:bijection})
making~(\ref{eqn:Irr}) commute.

From the constructions in~(\ref{sect:constructions})
(especially (\ref{subsect:G0})),
we obtain the following commutative diagram:
$$\begin{matrix}
  T_{F_\mu} &
  \rightarrow &
  T_0 &
  \leftarrow &
  {T_0}_{/ F_\lambda} &
  \xleftarrow{\overset{f_\lambda}{\cong}} &
  T_{F_\lambda} \\
 \hookdownarrow & &
 \hookdownarrow & &
 \hookdownarrow & &
 \hookdownarrow \\
  G_{F_\mu}  &
  \rightarrow &
  G_0 &
  \leftarrow &
  {G_0}_{/ F_\lambda} &
    &
  G_{F_\lambda} \\
 \downarrow & &
 \downarrow & &
 \downarrow & &
 \downarrow \\
  \Spec(F_\mu) &
  \rightarrow &
  \Spec(F) &
  \leftarrow &
  \Spec(F_\lambda) &
  = &
  \Spec(F_\lambda) 
\end{matrix}
$$
in which the left four squares
are cartesian.
Thanks to the fact that
the tori 
and the connected reductive groups
appearing above
are all split,
the functors of scalar extensions
(from $F$ to $F_\lambda$
and from $F$ to $F_\mu$ respectively)
yield
isomorphisms of character groups:
$$ \tX_{F_\mu}(T_{F_\mu})
   \xleftarrow{\cong}
   \tX_{F}(T_0)
   \xrightarrow{\cong}
   \tX_{F_\lambda}({T_0}_{/ F_\lambda})
$$
and bijections of sets of irreducible representations:
$$ \Irr_{F_\mu}(G_{F_\mu})
   \xleftarrow{\cong}
   \Irr_{F}(G_0)
   \xrightarrow{\cong}
   \Irr_{F_\lambda}({G_0}_{/ F_\lambda}),
$$
and these induce
the following commutative diagram
of Grothendieck rings:
$$\begin{matrix}
  \tK(T_{F_\mu}) &
  \xleftarrow{\cong} &
  \tK(T_0) &
  \xrightarrow{\cong} &
  \tK({T_0}_{/ F_\lambda}) &
  \xleftarrow{\overset{\tK(f_\lambda)}{\cong}} &
  \tK(T_{F_\lambda}) \\
 \scriptstyle{\ch_{G_{F_\mu}}} \ \hookuparrow & &
 \scriptstyle{\ch_{G_0}} \ \hookuparrow & &
 \scriptstyle{\ch_{{G_0}_{/ F_\lambda}}} \ \hookuparrow & &
 \hookuparrow \ \scriptstyle{\ch_{G_{F_\lambda}}} \\
  \tK(G_{F_\mu}) &
  \xleftarrow{\cong} &
  \tK(G_0) &
  \xrightarrow{\cong} &
  \tK({G_0}_{/ F_\lambda}) &
  &
  \tK(G_{F_\lambda}) \\
\end{matrix}
$$
We use the top row of isomorphisms
to identify
$\tK(T_{F_\mu})$
with
$\tK(T_{F_\lambda})$.
It follows that
to prove~(\ref{thm:bijection Irr}),
it suffices for us to construct
a map
\begin{equation*}
\tag{\ref{sect:reformulation}.1}
\label{eqn:Theta}
 \Theta:
   \Irr_{F_\mu}(G_{F_\mu})
   \rightarrow
   \Irr_{F_\lambda}(G_{F_\lambda})
\end{equation*}
which is bijective,
and which makes the diagram
\begin{equation*}
\tag{\ref{sect:reformulation}.2}
\label{eqn:Theta commute}
\begin{matrix}
  \tK(T_{F_\mu}) &
  = &
  \tK(T_{F_\lambda}) \\
 \scriptstyle{\ch_{G_{F_\mu}}} \ \hookuparrow & &
 \hookuparrow \scriptstyle{\ch_{G_{F_\lambda}}} \\
  \tK(G_{F_\mu}) &
  \xrightarrow{\Theta} &
  \tK(G_{F_\lambda})
\end{matrix}
\qquad\text{commute}.
\end{equation*}
Reformulating~(\ref{thm:bijection Irr}) in this way,
we see that
the two places $\lambda, \mu$
of the number field $F$
play symmetric roles.
In particular,
we see that
it suffices to construct
{\it some} map $\Theta$ in~(\ref{eqn:Theta})
which makes the diagram~(\ref{eqn:Theta commute})
commute;
it will automatically be bijective,
for if we reverse the roles of
$\lambda$ and $\mu$,
we would obtain
a map $\Theta'$
in the direction opposite to $\Theta$ in~(\ref{eqn:Theta}),
and the commutativity of the diagram~(\ref{eqn:Theta commute})
(and the analogous one for $\Theta'$)
would show that
$\Theta$ and $\Theta'$ are
inverses of each other.

\subsection{Constructions}
\label{sect:more constructions}
Recall from~(\ref{subsect:scalar extn})
and~(\ref{subsect:monodromy group}) that
$G_{F_\lambda}$ is the arithmetic mono\-dromy group
of the lisse $F_\lambda$-sheaf
$\Fcal_\lambda = \Lcal_\lambda\otimes_{E_\lambda} F_\lambda$.

\subsubsection{}
\label{subsect:monodromy rep}
The monodromy representation 
$[\Fcal_\lambda]$ of $\Fcal_\lambda$
factorizes as a composite homo\-morphism
$$ [\Fcal_\lambda] :
   \pi_1(X,\etabar)
   \xrightarrow{\omega_\lambda}
   G_{F_\lambda}
   \lhook\joinrel\xrightarrow{\rho_\lambda}
   \GL({\Fcal_\lambda}_\etabar),
$$
where
we let $\omega_\lambda$ denote
the continuous homomorphism from $\pi_1(X,\etabar)$
into $G_{F_\lambda}$ with a Zariski dense image,
and $\rho_\lambda$ is
the faithful tautological $F_\lambda$-rational representation
of $G_{F_\lambda}$
(cf.~(\ref{subsect:monodromy group})).

\subsubsection{}
\label{subsect:character gp}
From~(\ref{subsect:max torus}),
we see that
the maximal torus $T_{F_\lambda}$ of $G_{F_\lambda}$
is the Zariski closure of
the subgroup of $G_{F_\lambda}$
generated by $\omega_\lambda(\Frob_x)^\semisimp$,
where $x\in|X|$ is the point
chosen in~(\ref{subsect:closed point}).
The group of characters $\tX(T_{F_\lambda})$ of $T_{F_\lambda}$
is canonically isomorphic to
the subgroup $\Psi \subset F^\times$
generated by
the eigenvalues of $\Frob_x$
acting on ${\Fcal_\lambda}_\etabar$
(cf.~(\ref{sect:Frob tori})).
The Grothendieck ring $\tK(T_{F_\lambda})$ of $T_{F_\lambda}$
is therefore canonically isomorphic to
the group ring of $\Psi$ over $\bZ$.

\subsubsection{}
\label{subsect:Flambdatau}
Let $\tau: G_{F_\lambda} \rightarrow \GL_N$ be
an $F_\lambda$-rational representation of $G_{F_\lambda}$.
We may pull it back via $\omega_\lambda$
to obtain
a continuous $F_\lambda$-representation of $\pi_1(X,\etabar)$,
and hence
a lisse $F_\lambda$-sheaf on $X$,
which we will denote by $\Fcal_\lambda(\tau)$.

The restriction of $\tau$
to the maximal torus $T_{F_\lambda}$
decomposes as a direct sum of
characters of $T_{F_\lambda}$,
and this decomposition is encoded by
the element $\ch_{G_{F_\lambda}}(\tau)$
of $\tK(T_{F_\lambda})$.
Identifying $\tX(T_{F_\lambda})$ 
with the group $\Psi$ as in~(\ref{subsect:character gp}),
we see that
the characters of $T_{F_\lambda}$
appearing in this decomposition
are given by
the eigenvalues of $\Frob_x$
acting on
the lisse sheaf $\Fcal_\lambda(\tau)$;
in other words,
the element
$\ch_{G_{F_\lambda}}(\tau)\in\tK(T_{F_\lambda})$
is determined by
the ``inverse'' characteristic polynomial
$\det(1-T\,\Frob_x,\Fcal_\lambda(\tau))$
of $\Frob_x$ acting on
$\Fcal_\lambda(\tau)$.

\subsubsection{}
\label{subsect:absirred}
If $\tau$ is absolutely irreducible,
then so is the lisse $F_\lambda$-sheaf $\Fcal_\lambda(\tau)$;
this follows from the fact that
the monodromy representation of $\Fcal(\tau)$
is the composite of
the continuous homomorphism $\omega_\lambda$
(which has a Zariski dense image),
followed by
the representation $\tau$ of $G_{F_\lambda}$.

\subsubsection{}
\label{subsect:analog for mu}
The above discussion holds just as well
when $\lambda$ is replaced by $\mu$,
and we shall employ the analogous notation
for the objects constructed.
Note that by our hypotheses,
$\Fcal_\lambda$ and $\Fcal_\mu$
are $F$-compatible with each other;
likewise their duals
${\Fcal_\lambda}^\vee$ and ${\Fcal_\mu}^\vee$
are also $F$-compatible with each other.

The existence of the map $\Theta$ in~(\ref{eqn:Theta})
which makes the diagram~(\ref{eqn:Theta commute}) commute
is a consequence of the following:

\begin{theorem}
\label{prop:mu to lambda}
Given the constructions~(\ref{sect:constructions})
and~(\ref{sect:more constructions}) above,
for every
irreducible $F_\mu$-rational representation
$\theta_\mu\in\Irr_{F_\mu}(G_{F_\mu})$
of $G_{F_\mu}$,
there exists
an irreducible $F_\lambda$-rational representation
$\theta_\lambda\in\Irr_{F_\lambda}(G_{F_\lambda})$
of $G_{F_\lambda}$
such that
the $F_\mu$-sheaf $\Fcal_\mu(\theta_\mu)$
and
the $F_\lambda$-sheaf $\Fcal_\lambda(\theta_\lambda)$
are $F$-compatible.
Moreover,
one has
$$ \ch_{G_{F_\mu}}(\theta_\mu)
 = \ch_{G_{F_\lambda}}(\theta_\lambda)
 \qquad\text{equality in $\tK(T_{F_\mu}) = \tK(T_{F_\lambda})$}.
$$
\end{theorem}

\begin{proof}
We first work in the ``$\mu$-world''.
Recall that
the group $G_{F_\mu}$
is given with
the faithful tautological
$F_\mu$-rational representation $\rho_\mu$.
Let $\theta_\mu\in\Irr_{F_\mu}(G_{F_\mu})$ be
an irreducible $F_\mu$-rational representation
of $G_{F_\mu}$.
By a general result
in representation theory
(cf.~\cite{DemazureGabriel}, Chap.~II, \S2, Prop.~2.9),
there exist non-negative integers $a,b\in\bZ_{\ge 0}$
such that
$\theta_\mu$ occurs with positive multiplicity in
the tensor power representation
${\rho_\mu}^{\otimes a}
 \otimes
 {{\rho_\mu}^\vee}^{\otimes b}$.
Consider
the isotypic decomposition
$$ {\rho_\mu}^{\otimes a}
   \otimes
   {{\rho_\mu}^\vee}^{\otimes b}
   \cong
   \bigoplus_{\tau_\mu\in\Irr_{F_\mu}(G_{F_\mu})}
   {\tau_\mu}^{\oplus n(\tau_\mu)},
$$
where $n(\tau_\mu)$ is the multiplicity
(almost all of which are zero)
of $\tau_\mu$
in
${\rho_\mu}^{\otimes a}
 \otimes
 {{\rho_\mu}^\vee}^{\otimes b}$.
The isotypic decomposition of
the lisse $F_\mu$-sheaf
${\Fcal_\mu}^{\otimes a}
 \otimes
 {{\Fcal_\mu}^\vee}^{\otimes b}$
is then given by
\begin{equation*}
\tag{\ref{prop:mu to lambda}.1}
\label{eqn:ab decomp}
   {\Fcal_\mu}^{\otimes a}
   \otimes
   {{\Fcal_\mu}^\vee}^{\otimes b}
   \cong
   \bigoplus_{\tau_\mu\in\Irr_{F_\mu}(G_{F_\mu})}
   \Fcal_\mu(\tau_\mu)^{\oplus n(\tau_\mu)}.
\end{equation*}
Since both $\Fcal_\mu$ and ${\Fcal_\mu}^\vee$
are plain of characteristic~$p$,
so is
${\Fcal_\mu}^{\otimes a}
 \otimes
 {{\Fcal_\mu}^\vee}^{\otimes b}$.
Hence,
for each of the finitely many
$\tau_\mu\in\Irr_{F_\mu}(G_{F_\mu})$
with $n(\tau_\mu)>0$,
the lisse $F_\mu$-sheaf $\Fcal_\mu(\tau_\mu)$
is also
plain of characteristic~$p$.
Moreover,
because $\tau_\mu$ is
absolutely irreducible,
so is $\Fcal_\mu(\tau_\mu)$.
Therefore
each of these $\Fcal_\mu(\tau_\mu)$
satisfies the hypotheses of~(\ref{thm:indeplalg}).

We now ``pass from $\mu$ to $\lambda$''.
Applying~(\ref{thm:indeplalg}),
we see that
there is a finite extension $F'$ of $F$
and a place $\lambda$ of $F'$
lying over the given place $\lambda$ of $F$
such that
for each $\tau_\mu\in\Irr_{F_\mu}(G_{F_\mu})$
with $n(\tau_\mu)>0$,
there is an absolutely irreducible
lisse $F'_\lambda$-sheaf $\Hcal'(\tau_\mu)$ on $X$
which is $F$-compatible with
the lisse $F_\mu$-sheaf $\Fcal_\mu(\tau_\mu)$.
We can thus form
the lisse $F'_\lambda$-sheaf
\begin{equation*}
\tag{\ref{prop:mu to lambda}.2}
\label{eqn:Hcalprime}
   \Hcal' := 
   \bigoplus_{\tau_\mu\in\Irr_{F_\mu}(G_{F_\mu})}
   \Hcal'(\tau_\mu)^{\oplus n(\tau_\mu)}
\end{equation*}
(note that
the indexing set is still $\Irr_{F_\mu}(G_{F_\mu})$).
Comparing this with~(\ref{eqn:ab decomp}),
we see that
$\Hcal'$ is $F$-compatible with
the lisse $F_\mu$-sheaf
${\Fcal_\mu}^{\otimes a}
 \otimes
 {{\Fcal_\mu}^\vee}^{\otimes b}$.
As
$\Fcal_\mu$ (resp.~${\Fcal_\mu}^\vee$)
is $F$-compatible with
$\Fcal_\lambda$ (resp.~${\Fcal_\lambda}^\vee$),
the lisse $F'_\lambda$-sheaf $\Hcal'$
is also $F$-compatible with
the lisse $F_\lambda$-sheaf
${\Fcal_\lambda}^{\otimes a}
 \otimes
 {{\Fcal_\lambda}^\vee}^{\otimes b}$.
By \u{C}ebotarev's density theorem
and the trace comparison theorem of Bourbaki
(cf.~\cite{Bourbaki-AlgCh8}~\S12, no.~1, Prop.~3),
it follows that
$\Hcal'$ is isomorphic to
the $F'_\lambda$-scalar extension
of the lisse $F_\lambda$-sheaf
${\Fcal_\lambda}^{\otimes a}
 \otimes
 {{\Fcal_\lambda}^\vee}^{\otimes b}$.

From this observation,
we infer that
the monodromy representation of
the lisse $F'_\lambda$-sheaf $\Hcal'$
factors through the group
$(G_{F_\lambda})_{/ F'_\lambda}$
obtained from the group $G_{F_\lambda}$
by scalar extension to $F'_\lambda$.
Consequently,
each of the
absolutely irreducible
direct summands $\Hcal'(\tau_\mu)$ of $\Hcal'$
appearing in~(\ref{eqn:Hcalprime})
is a lisse $F'_\lambda$-sheaf
whose monodromy representation
also factors through the group
$(G_{F_\lambda})_{/ F'_\lambda}$.
In other words,
for each $\tau_\mu\in\Irr_{F_\mu}(G_{F_\mu})$
with $n(\tau_\mu)>0$,
there is
an $F'_\lambda$-rational representation $\tau_\lambda'$
of $(G_{F_\lambda})_{/ F'_\lambda}$
such that
the monodromy representation
of $\Hcal'(\tau_\mu)$
factorizes as
$$ [\Hcal'(\tau_\mu)] :
   \pi_1(X,\etabar)
   \xrightarrow{(\omega_\lambda)_{/ F'_\lambda}}
   (G_{F_\lambda})_{/ F'_\lambda}
   \xrightarrow{\tau_\lambda'}
   \GL(\Hcal'(\tau_\mu)_\etabar).
$$
The representation $\tau_\lambda'$
of $(G_{F_\lambda})_{/ F'_\lambda}$
is necessarily irreducible,
since this is so for $[\Hcal'(\tau_\mu)]$.

Now it remains for us
to descend
from $F'_\lambda$ to $F_\lambda$.
Since the group $G_{F_\lambda}$
is a split group over $F_\lambda$,
the representation $\tau_\lambda'$
is defined over $F_\lambda$:
there exists
an irreducible
$F_\lambda$-rational representation $\tau_\lambda$
of $G_{F_\lambda}$
whose scalar extension to $F'_\lambda$
is $\tau_\lambda'$.
Precomposing $\tau_\lambda$
with $\omega_\lambda$
yields
the continuous $F_\lambda$-representation
$\tau_\lambda\circ\omega_\lambda$
of $\pi_1(X,\etabar)$,
which is none other than
the monodromy representation of
the lisse $F_\lambda$-sheaf $\Fcal_\lambda(\tau_\lambda)$,
as constructed in~(\ref{subsect:Flambdatau}).
Since the $F'_\lambda$-scalar extension
of $\Fcal_\lambda(\tau_\lambda)$
is $\Hcal'(\tau_\mu)$,
we conclude that
$\Fcal_\lambda(\tau_\lambda)$
is $F$-compatible with
$\Fcal_\mu(\tau_\mu)$.

In summary,
to each
$\tau_\mu\in\Irr_{F_\mu}(G_{F_\mu})$
with $n(\tau_\mu)>0$,
we have constructed
a corresponding
$\tau_\lambda\in\Irr_{F_\lambda}(G_{F_\lambda})$
such that
$\Fcal_\lambda(\tau_\lambda)$
is $F$-compatible with
$\Fcal_\mu(\tau_\mu)$.
Since $a$ and $b$ are chosen so that
$n(\theta_\mu)>0$,
we may specialize to the case
of $\tau_\mu = \theta_\mu$,
and that proves
the first assertion of the theorem.
The second assertion follows
from the first
and the remarks in~(\ref{subsect:Flambdatau})
and~(\ref{subsect:absirred}).
This completes the proof
of~(\ref{prop:mu to lambda}),
and the proof of theorem~(\ref{thm:bijection Irr})
as well.
\end{proof}

\subsection{\it Proof of~(\ref{thm:main arith})(ii)}
\label{proof:main(ii)}
Fix $\lambda\in|F|_{\ne p}$
and use the isomorphism $f_\lambda$
given in~(\ref{cor:stronger main thm arith})
to identify the groups
$G_{F_\lambda}$ and ${G_0}_{/ F_\lambda}$.
We want to show that
under this identification,
the tautological representation $\rho_\lambda$
of $G_{F_\lambda}$
is isomorphic to
$\sigma_0\otimes_F F_\lambda$,
where $\sigma_0$ is
the $F$-rational representation of $G_0$
constructed in~(\ref{subsect:sigma0}).
It suffices to show that
the characters
$\ch_{G_{F_\lambda}}(\rho_\lambda)$
and
$\ch_{{G_0}_{/ F_\lambda}}(\sigma_0\otimes_F F_\lambda)$
are mapped to each other
under the isomorphism
$$ \tK(f_\lambda) : 
   \tK(T_{F_\lambda})
   \xrightarrow{\cong}
   \tK({T_0}_{/ F_\lambda}).
$$
From the construction
of $\sigma_0$ in~(\ref{subsect:sigma0})
and the argument in~(\ref{sect:reformulation}),
we see that
this is the same as showing
the equality of
$\ch_{G_{F_\lambda}}(\rho_\lambda)$
and
$\ch_{G_{F_\mu}}(\rho_\mu)$
as elements of
$\tK(T_{F_\lambda}) = \tK(T_{F_\mu})$.
In view of~(\ref{subsect:Flambdatau}),
this amounts to showing the equality
$$ \det(1-T\,\Frob_x,\Fcal_\lambda)
 = \det(1-T\,\Frob_x,\Fcal_\mu),
$$
but that follows from
the fact that
$\Fcal_\lambda$ and $\Fcal_\mu$
are $F$-compatible.
\qed




\providecommand{\bysame}{\leavevmode\hbox to3em{\hrulefill}\thinspace}

\end{document}